\documentclass[onecolumn]{autart}    % Enable this line and disable the 
\usepackage{graphicx}          % Include this line if your 
\usepackage{amsmath}
\usepackage{amsfonts}
\usepackage{amssymb}
\usepackage{float}
\usepackage{graphicx}
\usepackage{mathtools}
\usepackage{booktabs}
\usepackage{xfrac}
\theoremstyle{definition}
\newtheorem{theorem}{Theorem}

\theoremstyle{definition}
\newtheorem{definition}{Definition}[section]
\theoremstyle{definition}
\newtheorem{assumption}{Assumption}
\newtheorem{proposition}{Proposition}
\newtheorem{lemma}{Lemma}
\newtheorem{example}{Example}

\usepackage{pgfplots}
\pgfplotsset{compat=newest}
\usetikzlibrary{plotmarks}
\usetikzlibrary{arrows.meta}
\usepgfplotslibrary{patchplots}
\usepackage{grffile}

\newcommand{\R}{{\mathbb R}}
\DeclareMathOperator{\Ima}{Im}

\newcommand{\added}{\textcolor{black}}
\newcommand{\addedTwo}{\textcolor{black}}

\usepackage{xcolor}

\begin{document}

\begin{frontmatter}

\title{Approximate abstractions of control systems with an application to aggregation} % Title, preferably not more 

\thanks[footnoteinfo]{This work was supported in part by a NDSEG Graduate Fellowship, NSF grants  CNS-1446145 and CNS-1545116, the AFOSR grant FA9550-18-1-0253, the H2020 ERC Starting Grant AutoCPS (grant agreement No. 804639), and the German Research Foundation (DFG) through the grant ZA 873/1-1. The material in this paper was partially presented at the 2018 American Control Conference, June 27 - 29, Milwaukee, WI, USA. Corresponding author S. W. Smith. \textsuperscript{\textcopyright} 2020 the authors. Licensed under a Creative Commons Attribution-NonCommercial-NoDerivatives 4.0 International License. \\
See https://creativecommons.org/licenses/by-nc-nd/4.0/ for further information.}

\author[Berkeley]{Stanley W. Smith}\ead{swsmth@berkeley.edu},    % Add the 
\author[Berkeley]{Murat Arcak}\ead{arcak@berkeley.edu},               % e-mail address 
\author[Boulder,Munich]{Majid Zamani}\ead{majid.zamani@colorado.edu}  % (ead) as shown

\address[Berkeley]{Department of Electrical Engineering and Computer Sciences, University of California, Berkeley}  % Please supply                                              
\address[Boulder]{Department of Computer Science, University of Colorado Boulder, USA}

\address[Munich]{Department of Computer Science, Ludwig Maximilian University of Munich, Germany}             % full addresses

\begin{keyword}                           % Five to ten keywords,  
Approximate Abstractions, Practical Simulation/Storage Functions, Aggregation             % chosen from the IFAC 
\end{keyword}                             % keyword list or with the 
                                          % help of the Automatica 
                                          % keyword wizard

\begin{abstract} Previous approaches to constructing abstractions for control systems rely on geometric conditions or, in the case of an interconnected control system, a condition on the interconnection topology. Since these conditions are not always satisfiable, we relax the restrictions on the choice of abstractions, instead opting to select ones which nearly satisfy such conditions via optimization-based approaches. To quantify the resulting effect on the error between the abstraction and concrete control system, we introduce the notions of practical simulation functions and practical storage functions. We show that our approach facilitates the procedure of aggregation, where one creates an abstraction by partitioning agents into aggregate areas. We demonstrate the results on an application where we regulate the temperature in three separate zones of a building.
\end{abstract}

\end{frontmatter}

\section{Introduction}
The synthesis of controllers for dynamical systems enforcing complex logic properties, e.g. those expressed as linear or signal temporal logic (LTL/STL) formulas \cite{katoen,donze2013signal}, is hampered by computational challenges. One way of tackling the design complexity is by employing abstractions, which are simpler representations of original systems with the property that controllers designed for them to enforce desired properties can be refined to the ones for the concrete systems. The errors suffered in this controller synthesis detour can be quantified a priori. The abstraction is called finite if its set of states is finite, and infinite otherwise. In this paper, we only deal with infinite abstractions.

Abstractions of non-stochastic dynamical systems has a long history. Examples of such results include constructive procedures for the construction of infinite abstractions of linear control systems using exact simulation relations \cite{vdS04}. In contrast to the exact notions, the results in \cite{GP09} provide an approach for the construction of infinite abstractions of linear control systems using approximate simulation relations based on \textit{simulation functions}. The construction schemes proposed in \cite{vdS04,GP09} are monolithic in the sense that infinite abstractions are constructed from the complete system model. Compositional construction of approximate abstractions for the interconnection of two subsystems is studied in \cite{girard2013composition} using small gain type conditions. This result was extended in \cite{RZ1} to networks of systems, again with small gain type reasoning. The recent result in \cite{ZA3} employs broader dissipativity methods for constructing approximate abstractions for networks. 

The infinite abstractions discussed here are also related to the rich theory of model order reduction, which seeks abstractions with reduced state-space dimensions \cite{Antoulas05,SM09}. However, the model mismatch in \cite{Antoulas05,SM09} is established with respect to $\mathcal{H}_2/\mathcal{H}_\infty$ norms whereas we use notions of simulation functions to derive $\mathcal{L}_\infty$ error bounds, which are crucial to reason about complex logic properties, e.g., LTL or STL formulas \cite{katoen,donze2013signal}.    

The aforementioned results on the construction of exact or approximate infinite abstractions, \cite{vdS04,GP09,RZ1,ZA3}, require restrictive geometric conditions which, in some cases, are satisfied only when the state dimensions of the abstraction and the original system are the same (i.e., no order reduction).

In this work, we address this shortcoming as follows. We first show that, when constructing an abstraction monolithically, one can relax the geometric conditions appearing in \cite{vdS04,GP09,RZ1,ZA3}. We quantify the effect of this relaxation via a nonnegative function which can be bounded in a formal synthesis of the abstract controller. To translate this bound into one on the error between the concrete system and its abstraction, we modify the definition of simulation functions from \cite{GP09} to that of \textit{practical simulation functions}, which include the nonnegative function in the upper bound on their derivative.

Next, we show that when constructing an abstraction in a compositional manner, one can also relax a restrictive condition on the interconnection topology from \cite{RZ1,ZA3}. We show that this relaxation greatly expands the domain of applicability of model order reduction via \textit{aggregation}, where one creates an abstraction by partitioning agents into aggregate areas. In addition, our construction utilizes a modified version of storage functions from \cite{ZA3}, which we refer to as \textit{practical storage functions}. This notion allows us to accommodate heterogeneity in the agent models in aggregation.

\added{The flexibility of our approach greatly broadens the applicability of infinite abstractions, including their usage in formal control synthesis procedures. Indeed, it was previously difficult and at times intractable to find an infinite abstraction satisfying the aforementioned geometric conditions. Thus, our method overcomes a significant limitation of abstraction-based controller design by allowing one to instead use an approximate abstraction which need not satisfy such conditions. The additional error introduced by this approach can then be quantified with our newly introduced notion of a practical simulation function.}

The paper is organized as follows. In Section \ref{controlSys}, we introduce the class of control systems and corresponding abstractions studied in the paper. We show in Section \ref{monolithic} how one can construct an abstraction in a monolithic manner for the class of linear systems. The discussion in Section \ref{monolithic} is based on the preliminary work in \cite{smith2018hierarchical}; however, the content after Section \ref{monolithic} is entirely new. In Section \ref{Composition-Section}, we consider a class of interconnected control systems, and present a result on the compositional construction of an abstraction for such systems. In Section \ref{aggregationSection}, we show how our theory can aid in the procedure of aggregation, and include an example in building temperature regulation in Section \ref{temperatureExample}. We conclude with final remarks in Section \ref{conclusionSection}. All proofs are given in the Appendix.

\section{Control Systems} \label{controlSys}

\subsection{Notation.}
We denote the set of real numbers as $\R$, and write the set of positive and nonnegative real numbers as $\R_{>0}$ and $\R_{\geq 0}$, respectively. For $a, b \in \R$ with $a \leq b$, we denote with $(a,b)$ the open interval from $a$ to $b$. The $n$-dimensional Euclidean space is denoted with $\R^n$. We use $\mathbf{1}_n$ and $\mathbf{0}_n$ to denote the $n$-dimensional vector with all entries equal to $1$ and $0$, respectively. The vector space of matrices with $n$ rows and $m$ columns is represented by $\R^{n \times m}$. We use $I_n$ to denote the identity matrix with $n$ rows and columns. The concatenation of vectors $x_i \in \R^{n_i}$ for $i = 1, \dots, N$ is given by $[x_1; x_2; \dots; x_N] \in \R^n$, where $n = \sum_{i=1}^N n_i$. Similarly, the block-diagonal concatenation of matrices $P_i \in \R^{m_i \times n_i}$ for $i = 1, \dots, N$ is written as $\text{diag}(P_1, \dots, P_N) \in \R^{m \times n}$, where $m$ and $n$ are defined in the same way. The null space of a matrix $P \in \R^{m \times n}$ is given by $\mathcal{N}(P) := \{x \in \R^{n} : Px = \mathbf{0}_m\}$. Furthermore, $\|P\|_F$ and $\mathbf{tr}(P)$ refer to the Frobenius norm and trace of $P$, respectively. The map $\|\cdot\| : \R^{n \times m} \to \R_{\geq 0}$ refers to the Euclidean norm when the argument is a vector, and the matrix norm induced by the Euclidean norm when the argument is a matrix. For a symmetric matrix $P \in \R^{n \times n}$, we use $\lambda_{\text{min}}(P)$ and $\lambda_{\text{max}}(P)$ to denote the minimum and maximum eigenvalue of $P$, respectively. We denote the Kronecker product of matrices $A \in \R^{m \times n}$ and $B \in \R^{p \times q}$ as $A \otimes B \in \R^{mp \times nq}$.

A continuous function $\alpha : \R_{\geq 0} \to \R_{\geq 0}$ belongs to class $\mathcal{K}$ if it is strictly increasing and $\alpha(0) = 0$; furthermore, $\alpha : \R_{\geq 0} \to \R_{\geq 0}$ belongs to class $\mathcal{K}_{\infty}$ if $\alpha \in \mathcal{K}$ and $\alpha(s) \to \infty$ as $s \to \infty$. A continuous function $\beta : \R_{\geq 0} \times \R_{\geq 0} \to \R_{\geq 0}$ is said to belong to class $\mathcal{K} \mathcal{L}$ if, for each fixed $t$, the map $\beta(r, t)$ belongs to class $\mathcal{K}$ with respect to $r$ and, for each fixed nonzero $r$, the map $\beta(r,t)$ is decreasing with respect to $t$ and $\beta(r,t) \to 0$ as $t \to \infty$. Lastly, for a measurable function $f : \R_{\geq 0} \to \R^{n}$, we use $\|f\|_{\infty}$ to indicate $\sup_{t \geq 0} \|f(t)\|$.

\subsection{Control systems and their abstractions.}
We first define the class of control systems studied in this paper:

\begin{definition}
A control system $\Sigma$ is a tuple $\Sigma = (\R^n, \R^m, f, \R^q, h)$, where $\R^{n}$, $\R^m$, and $\R^q$ are the state, input, and output spaces, respectively. The evolution of the state and output trajectories are governed by
\begin{align}
\Sigma : 
\begin{cases}
\dot{\xi}(t) = f(\xi(t), \upsilon(t)), \\
\zeta(t) = h(\xi(t)),
\end{cases} \nonumber %\label{full}
\end{align}
where $f: \R^n \times \R^m \to \R^n$ is locally Lipschitz, and we refer to $h: \R^n \to \R^q$ as the output map.
\end{definition}
We denote by $\xi_{x\upsilon}(t)$ the state reached at time $t$
under the input $\upsilon:\R_{\geq0}\to\R^m$ from the initial condition $x=\xi_{x\upsilon}(0)$; the state $\xi_{x\upsilon}(t)$ is
uniquely determined due to the assumptions on $f$ \cite{sontag1}. We also denote by $\zeta_{x\upsilon}(t)$ the corresponding output value of $\xi_{x\upsilon}(t)$, i.e. $\zeta_{x\upsilon}(t)=h(\xi_{x\upsilon}(t))$.

When the dimension of the state space is large, one can avoid the computational burden of a direct controller synthesis for $\Sigma$ by introducing an abstraction $\hat{\Sigma}$, potentially with a smaller state-space dimension $\hat{n}$. Typically, the abstraction $\hat{\Sigma}$ is related to the concrete system $\Sigma$ via a \textit{simulation function} \cite{GP09}, which enables one to bound the error between the outputs of the two systems. We now define a modified version of simulation functions, which we refer to as {\it practical} simulation functions:

\begin{definition} \label{simFunDef}
Consider a control system $\Sigma = (\R^n, \R^m, f, \R^q, h)$ with corresponding abstraction $\hat{\Sigma} = (\R^{\hat{n}}, \R^{\hat{m}}, \hat{f}, \R^q, \hat{h})$. Let $V: \mathbb{R}^n \times \mathbb{R}^{\hat{n}} \to \mathbb{R}_{\geq 0}$ be a continuously differentiable function and $v : \R^n \times \R^{\hat{n}} \times \R^{\hat{m}} \to \R^m$ a locally Lipschitz function. We say that $V$ is a practical simulation function from $\hat{\Sigma}$ to $\Sigma$ with an associated interface $v$ if there exist $\nu, \eta \in \mathcal{K}_\infty$, $\rho \in \mathcal{K} \cup \{0\}$, and $\Delta : \R^{\hat{n}} \to \R_{\geq 0}$ such that for all $x$, $\hat{x}$, and $\hat{u}$ we have
\begin{equation} \label{output}
\nu(\|h(x) - \hat{h}(\hat{x})\|) \ \leq V(x, \hat{x})
\end{equation}
and
\begin{align}
& \frac{\partial V(x, \hat{x})}{\partial x} f(x, v(x, \hat{x}, \hat{u})) + \frac{\partial V(x, \hat{x})}{\partial \hat{x}} \hat{f}(\hat{x}, \hat{u}) \leq - \eta(V(x, \hat{x})) + \rho(\|\hat{u}\|) + \Delta(\hat{x}). \label{iss}
\end{align}
\end{definition}
Here, we modified the definition of simulation functions to include a nonnegative term $\Delta(\hat{x})$ in the upper bound of their derivatives. \added{Thus, when $\Delta(\hat{x}) = 0$ we refer to $V(x, \hat{x})$ as a simulation function. We note that the associated interface $v(x, \hat{x}, \hat{u})$ helps to achieve \eqref{iss} and, in particular, can be used to reduce the term $\Delta(\hat{x})$ as much as possible}. The usefulness of $\Delta(\hat{x})$ will become apparent in Section \ref{monolithic}, where we show that its addition allows one to relax the geometric conditions typically required in the construction of infinite abstractions. \added{To further motivate the addition of the term $\Delta(\hat{x})$, we provide an example of a system and abstraction which admit a practical simulation function as in Definition 2.2:
\begin{example}
Consider the control system
\begin{equation*}
\Sigma: \begin{cases}
\begin{aligned}
\begin{pmatrix}
\dot{\xi}_1(t) \\
\dot{\xi}_2(t)
\end{pmatrix} &=
\begin{pmatrix}
-1.5 \xi_1^3 (t) + \upsilon(t) \\
-\xi_2^3(t) + \upsilon(t) 
\end{pmatrix}, \\[5pt]
\zeta(t) &= (\xi_1 (t) + \xi_2 (t))/2,
\end{aligned}
\end{cases}
\end{equation*}
with $\xi_1(t), \xi_2(t), \zeta(t) \in \R$, and where $\xi_1(t)$ and $\xi_2(t)$ are aggregated into a single state variable $\hat{\xi}(t) \in \R$ governed by
\begin{equation*}
\hat{\Sigma}: \begin{cases}
\dot{\hat{\xi}}(t) = - 1.5{\hat{\xi}}^3 (t) + \hat{\upsilon}(t), \\
\hat{\zeta}(t) = \hat{\xi}(t),
\end{cases}
\end{equation*}
with $\hat{\zeta}(t) \in \R$. Then, by defining the associated interface
\begin{equation*}
v(x, \hat{x}, \hat{u}) = \hat{u},
\end{equation*}
we have that $V(x, \hat{x}) := (1/2) (x - \hat{x} \mathbf{1}_2)^T (x - \hat{x} \mathbf{1}_2)$ is a practical simulation function from $\hat{\Sigma}$ to $\Sigma$ since one can verify that
\begin{align*}
((x_1 + x_2)/2 - \hat{x})^2 \leq V(x, \hat{x})
\end{align*}
and
\begin{align*}
\dot{V}(x, \hat{x}) &= (x - \hat{x} \mathbf{1}_2)^T (\dot{x} - \dot{\hat{x}} \mathbf{1}_2) \\
&= \begin{bmatrix}
x_1 - \hat{x} \\
x_2 - \hat{x}
\end{bmatrix}^T \begin{bmatrix}
-1.5 x_1^3 + \hat{u} - (-1.5\hat{x}^3 + \hat{u}) \\
- x_2^3 + \hat{u} - (-1.5\hat{x}^3 + \hat{u})
\end{bmatrix} \\
% &= \begin{bmatrix}
% x_1 - \hat{x} \\
% x_2 - \hat{x}
% \end{bmatrix}^T \begin{bmatrix}
% -(x_1^3 - \hat{x}^3) - x_1^3 \\
% -(x_2^3 + \hat{x}^3) 
% \end{bmatrix} \\ 
% & \leq -2V^2(x, \hat{x}) + (x_1 - \hat{x}) \cdot x_1^3 \\
% & \leq -2V^2(x, \hat{x}) + (1/4) (x_1 - \hat{x})^4 + (3/4) \hat{x}^4 \\ 
& \leq -\frac{1}{8} V^2(x, \hat{x}) + \frac{3}{8} \hat{x}^4
\end{align*}
hold. Thus, we have that \eqref{output} and \eqref{iss} from Definition 2.2 are satisfied with $\nu(s) := s^2, \ \eta(s) := (1/8) s^2, \ \rho(s) = 0$, and $\Delta(\hat{x}) := (3/8) \hat{x}^4$.
\end{example}
}

The next theorem shows the usefulness of a practical simulation function by providing a bound on the error between the output behaviors of control systems to those of their
abstractions.

\begin{theorem} \label{simulationFunctionTheorem}
Consider a system $\Sigma = (\R^n,\R^m,f,\R^q, h)$ with corresponding abstraction $\hat{\Sigma} = (\R^{\hat{n}},\R^{\hat{m}},\hat{f},\R^q,\hat{h})$, and let $V$ be a practical simulation function from $\hat{\Sigma}$ to $\Sigma$. Then, there exists a class $\mathcal{K} \mathcal{L}$ function $\beta$ and class $\mathcal{K}$ functions $\gamma_1, \gamma_2$ such that for any measurable $\hat{\upsilon} : \R_{\geq 0} \to \R^{\hat{m}}$ and $x \in \R^n$, $\hat{x} \in \R^{\hat{n}}$, there exists a measurable $\upsilon : \R_{\geq 0} \to \R^m$ via the associated interface $v$ such that the following bound holds for all $t \in \R_{\geq 0}$:
\begin{align*}
\| \zeta_{x\upsilon}(t) - \hat{\zeta}_{\hat x\hat\upsilon}(t) \| \leq \beta(V(x,\hat{x}), t) + \gamma_1( \| \hat{\upsilon} \|_{\infty} ) + \gamma_2( \| \Delta(\hat{\xi}_{\hat x\hat\upsilon}) \|_{\infty} ). \nonumber
\end{align*}
\end{theorem}
The proof of Theorem \ref{simulationFunctionTheorem} is similar to the one of Theorem 3.5 in \cite{ZA3} and is omitted here due to lack of space.

%\section{Abstraction Synthesis for a Class of Control Systems} 
\section{Abstraction Synthesis for Linear Systems}
\label{monolithic}

%\subsection{Linear Systems} \label{linear}
%In this section we briefly recall a result from \cite{smith2018hierarchical} on monolithic abstraction synthesis for the class of affine control systems. For the sake of clarity, 
To demonstrate the relaxation of geometric constraints, here we adapt our approach to 
%the class of 
linear control systems 
%$\Sigma = (\R^n, \R^m, f, \R^q, h)$, with functions $f$ and $h$ described by
\begin{align} \label{linearConcrete}
\Sigma :
\begin{cases}
\dot{\xi}(t) = A\xi(t) + B\upsilon(t), \\
\zeta(t) = C\xi(t),
\end{cases}
\end{align}
where $A \in \R^{n \times n}$, $B \in \R^{n \times m}$, $C \in \R^{q \times n}$, and the pair $(A,B)$ is stabilizable. Our goal is to represent \eqref{linearConcrete} with an abstract control system 
%$\hat\Sigma = (\R^{\hat n}, \R^{\hat m}, \hat f, \R^q, \hat h)$ with functions $\hat f$ an $\hat h$ given as
\begin{align} \label{linearAbstract}
\hat{\Sigma} :
\begin{cases}
\dot{\hat{\xi}}(t) = \hat{A} \hat{\xi}(t) + \hat{B} \hat{\upsilon}(t), \\
\hat{\zeta}(t) = \hat{C} \hat{\xi}(t),
\end{cases}
\end{align}
where $\hat{A} \in \R^{\hat{n} \times \hat{n}}$, $\hat{B} \in \R^{\hat{n} \times \hat{m}}$, and $\hat{C} \in \R^{q \times \hat{n}}$. It has been shown in \cite[Theorem 2]{GP09} that if one can find matrices $P \in \R^{n \times \hat{n}}$ and $Q \in \R^{m \times \hat{n}}$ such that $\hat{C} = C P$, and the condition
\begin{equation} \label{linearCondition}
AP = P \hat{A} - BQ
\end{equation}
holds, then there exists a \added{practical} simulation function from $\hat{\Sigma}$ to $\Sigma$ with an associated interface given by
\begin{equation} \label{linearInterface}
v(x, \hat{x}, \hat{u}) = K (x - P \hat{x}) + Q \hat{x} + R \hat{u}
\end{equation}
where the matrix $K \in \R^{m \times n}$ in \eqref{linearInterface} is a feedback gain to be designed and $R \in \R^{m \times \hat{m}}$ is selected to minimize $\|BR - P\hat{B}\|$. As alluded to previously, the requirement \eqref{linearCondition} can be restrictive in general. Indeed, the following lemma, quoted from \cite[Lemma 2]{GP09}, provides the geometric conditions on $P$ such that \eqref{linearCondition} is satisfiable:
\begin{lemma}  
For given matrices $A$, $P$, and $B$, there exist matrices $\hat{A}$ and $Q$ satisfying \eqref{linearCondition} if and only if
\begin{align}
\Ima(A P) \subseteq \Ima(P) + \Ima(B). \nonumber %\label{geometric}
\end{align}
\end{lemma}
%The proof is straightforward and is omitted. 
To address the restriction implicit in \eqref{linearCondition}, we propose a relaxation by allowing a nonzero residual term given by
\begin{equation}
D := AP - P\hat{A} + BQ. \nonumber
\end{equation}
The effect of a nonzero matrix $D$ is seen by examining the dynamics of the error $e(t) := \xi (t) - P \hat{\xi}(t)$, which become
\begin{equation} \label{linearErrorDyn}
\dot{e}(t) = (A + BK) e(t) + D \hat{\xi}(t) + (BR - P\hat{B}) \hat{\upsilon}(t)
\end{equation}
where
\begin{equation} \label{linearDisturbance}
D \hat{\xi}(t) + (BR - P\hat{B}) \hat{\upsilon}(t)
\end{equation}
is treated as a disturbance. Thus, by relaxing \eqref{linearCondition}, we have introduced a new term depending on $\hat{\xi}$ into the disturbance \eqref{linearDisturbance}, which previously only depended on $\hat{\upsilon}$. 

%In the next section, we design the feedback gain $K$ to mitigate the effect of this disturbance.

%\subsection{Robustness Analysis} \label{robust}
%In this section, we provide robustness bounds on $e$. We reduce 

We next design the feedback gain $K$ to mitigate the effect of this disturbance. To this end we rewrite \eqref{linearErrorDyn} as
\begin{equation} \label{efinal}
\dot{e}(t) = (A + BK) e(t) + W d(t)
\end{equation}
where we have defined
\begin{equation} \label{dform}
W := \begin{bmatrix}
I & BR - P \hat{B}
\end{bmatrix}, \quad d := \begin{bmatrix}
D \hat{\xi} \\ \hat{\upsilon}
\end{bmatrix},
\end{equation}
where $I$ is the identity matrix of appropriate size. 
The magnitude of $d$ can be bounded by placing constraints on $D \hat{\xi}$ and $\hat{\upsilon}$, to be respected for all $t \geq 0$. This can be done by introducing an appropriate STL specification for $\hat{\Sigma}$ which constrains $D \hat{\xi}$ and $\hat{\upsilon}$, and then synthesizing a control law $\hat{\upsilon}$ such that the resulting trajectories of $\hat{\Sigma}$ satisfy said specification - known as a formal synthesis procedure. In this paper, we apply a formal synthesis procedure utilizing \textit{model predictive control} (MPC) \cite{raman2017}; MPC is well known for being able to handle such constraints. Note that we do not need to constrain $\hat{\xi}$ itself to be small, but rather the value of $D \hat{\xi}$. For example, in a motion coordination application in \cite{smith2018hierarchical}, $D \hat{\xi}$ yields relative positions and the constraints do not unreasonably restrict the absolute positions contained in the vector $\hat{\xi}$. %\Majid{"formal synthesis procedure" is not very precise in the previous sentences. We never talked about it before! Please elaborate more on the last two sentences.}

\addedTwo{We remark that using MPC to design $\hat{\upsilon}$ requires discretization of the dynamics \eqref{linearAbstract}. This is important to note, in particular, since this implies the constraints on $d$ in \eqref{dform} will only hold at each sampling instant. Thus, we must establish a growth bound on each component of $d$ in order to characterize its inter-sample behavior. For $\hat{\xi}$, one can impose constraints such that $\hat{A} \hat{\xi}$ and $\hat{B} \hat{\upsilon}$ are bounded, and then subsequently bound $\dot{\hat{\xi}}$ from \eqref{linearAbstract}. Furthermore, since $\hat{\upsilon}$ is a zero-order hold signal its derivative between samples is zero. Combining these facts to provide a bound on $d$, we ensure the quality of the abstraction $\hat{\Sigma}$.}

\addedTwo{After designing $\hat{\upsilon}$, our goal} becomes to design $K$ to minimize the $\mathcal{L}_\infty$ gain from $d$ to error $e$. Since \eqref{efinal} is linear, an estimate for this gain is obtained by finding a bound $\overline{e} := \|e\|_\infty$ when $\overline{d} := \|d\|_\infty \leq 1$. We pursue this by numerically searching for $U = U^T > 0$ such that the ellipsoid $\mathcal{E} = \{e : e^T U e \leq 1 \}$ is invariant. This results in $\overline{e} = 1 / \sqrt{\lambda_{\min}(U)}$, since this is the radius of the smallest ball enclosing $\mathcal{E}$. The following optimization problem combines the search for $U$ with a simultaneous search for a $K$ that minimizes $\overline{e}$. Its derivation is similar to Section 6.1.3 of \cite{boyd1994linear} and is omitted here due to lack of space. \\[5pt]
\textbf{Optimization Problem 1:}
\begin{align}
\text{minimize } & \quad \beta \quad \text{over} \ Z := U^{-1}, \ Y := KZ, \nonumber \\[5pt]
\text {subject to } & \quad Z \leq \beta I, \label{constr1} \\[5pt]
& \quad X(Z, Y, \alpha) \leq 0, \label{lmiconstr1}
\end{align}
where
\begin{align*}
X(Z, Y, \alpha) :=
\begin{bmatrix}
A Z + Z A^T + Y^T B^T + B Y + \alpha Z & W \\
W^T & -\alpha I
\end{bmatrix},
\end{align*}
which is an LMI in $Z$ and $Y$ if the scalar $\alpha > 0$ is fixed. \addedTwo{In particular, by minimizing $\beta$ and imposing \eqref{constr1}, we are effectively maximizing $\lambda_{\min}(U)$. Here, this is equivalent to minimizing the error bound since $\overline{e} = 1 / \sqrt{\lambda_{\min}(U)}$}. The next theorem states that a solution to Optimization Problem 1 yields a practical simulation function from $\hat{\Sigma}$ to $\Sigma$.

\begin{theorem} \label{linearTheorem}
Suppose that $U$ and $K$ are a solution to Optimization Problem 1, and $\hat{C}$ in \eqref{linearAbstract} satisfies $\hat{C} = CP$. Then $V(x, \hat{x}) := (x - P \hat{x})^T U (x - P \hat{x})$ is a practical simulation function from $\hat{\Sigma}$ to $\Sigma$ with an associated interface $v(x, \hat{x}, \hat{u})$ as in \eqref{linearInterface}.
\end{theorem}
% \Majid{Double check the proof in the appendix and fix the typos! There, you are using $M$ whereas in the theorem you are using $U$.}

\addedTwo{As mentioned in Theorem \ref{simulationFunctionTheorem}, the practical simulation function $V(x, \hat{x})$ bounds the error between the outputs of $\Sigma$ and $\hat{\Sigma}$. This allows us to translate guarantees on $\hat{\Sigma}$ to weakened guarantees on $\Sigma$. For example, if one designs a controller enforcing a set $\hat \Omega$ to be invariant for $\hat \Sigma$, then the refined controller makes $\Omega^{\bar e}$ invariant for $\Sigma$, where in this case $\Omega^{\bar e} := \{e + P \hat{x} : \|e\|_\infty \leq \bar{e}, \ \hat{x} \in \hat{\Omega} \}$. The question then becomes how to obtain a small bound on $e$ so that the desired behavior is realized on $\Sigma$. A rigorous procedure for doing so is not the main focus of this paper, but is explored in \cite{yin2019optimization}. Here, we simply focus on improving the error bound via the two steps outlined in this section: first, by designing $\hat{\upsilon}$ to restrict $d$, and second, by using the interface $v(x, \hat{x}, \hat{u})$ to reduce the gain from $d$ to $e$. Our procedure is oriented towards control synthesis, as our goal is to move from designing an abstract controller towards designing a concrete one. In verification, where one wants to verify behavior correctness via abstraction, these steps cannot be applied in the reverse direction to reduce error, which could result in poor abstraction quality. Thus, we remark that our approach cannot be extended to verification in a straightforward way.}

\section{Compositionality} \label{Composition-Section}

\subsection{Interconnected control systems}

In this section we propose an approach to construct an abstraction and corresponding practical simulation function for a class of interconnected control systems. In particular, we show how to do so 
%in a compositional manner by leveraging knowledge of 
by composing the abstractions 
%and \textit{storage functions} for 
of the subsystems.
%comprising the interconnected control system. 
We start by  defining the class of subsystems that we consider:

\begin{definition} \label{controlSubsystem}
A control subsystem $\Sigma$ is a tuple $\Sigma = (\R^n, \R^m, \R^p, f, \R^{q_1}, \R^{q_2}, h_1, h_2)$, where $\R^n$, $\R^m$, $\R^p$, $\R^{q_1}$, and $\R^{q_2}$ are the state, external input, internal input, external output, and internal output spaces, respectively. The evolution of the state and output trajectories are governed by the equations
\begin{equation}
\Sigma : \begin{cases}
\begin{aligned}
\dot{\xi}(t) &= f(\xi(t), \upsilon(t), \omega(t)), \\
\zeta_{1}(t) &= h_{1}(\xi(t)), \\
\zeta_{2}(t) &= h_{2}(\xi(t)), 
\end{aligned}
\end{cases} \nonumber
\end{equation}
where $f: \R^n \times \R^m \times \R^p \to \R^n$ and $h_2 : \R^n \to \R^{q_2}$ are locally Lipschitz.  We refer to $h_1 : \R^n \to \R^{q_1}$ and $h_2 : \R^n \to \R^{q_2}$ as the \textit{external} and \textit{internal} output maps, respectively.
\end{definition}

Similar to a practical simulation function, a \textit{storage function} \cite{ZA3} can be used to relate a control subsystem $\Sigma$ to its abstraction $\hat{\Sigma}$ by describing a dissipativity property of the error dynamics.

\begin{definition} \label{storageFnDef}
Consider a control system $\Sigma = (\R^n, \R^m, \R^p, f, \R^{q_1}, \R^{q_2}, h_1, h_2)$ and corresponding abstraction $\hat{\Sigma} = (\R^{\hat{n}}, \R^{\hat{m}}, \R^{\hat{p}}, \hat{f}, \R^{q_1}, \R^{\hat{q}_2}, \hat{h}_1, \hat{h}_2)$. Let $V: \mathbb{R}^n \times \mathbb{R}^{\hat{n}} \to \mathbb{R}_{\geq 0}$ be a continuously differentiable function and $v : \R^n \times \R^{\hat{n}} \times \R^{\hat{m}} \to \R^m$ a locally Lipschitz function. We say that $V$ is a practical storage function from $\hat{\Sigma}$ to $\Sigma$ if there exist $\nu, \eta \in \mathcal{K}_\infty$, $\rho \in \mathcal{K} \ \cup \ \{ 0 \}$, a function $\Delta: \R^{\hat{n}} \to \R_{\geq 0}$, matrices $W$, $\hat{W}$, $H$ of appropriate dimensions, and matrix $X = X^T$ of appropriate dimension with conformal block partitions $X^{11}$, $X^{12}$, $X^{21}$, and $X^{22}$, such that for any $x \in \mathbb{R}^n$, $\hat{x} \in \mathbb{R}^{\hat{n}}$, $\hat{u} \in \mathbb{R}^{\hat{m}}$, $\hat{w} \in \mathbb{R}^{\hat{p}}$, and $w \in \mathbb{R}^p$ we have
\begin{equation} \label{storageBound}
\nu( \| h_1(x) - \hat{h}_1(\hat{x}) \| ) \ \leq V(x, \hat{x}) \nonumber
\end{equation}
and 
\begin{align}
\frac{\partial V(x, \hat{x})}{\partial x} f&(x, v(x, \hat{x}, \hat{u}), w) + \frac{\partial V(x, \hat{x})}{\partial \hat{x}} \hat{f}(\hat{x}, \hat{u}, \hat{w}) \nonumber \\
& \leq -\eta(V(x, \hat{x})) + \rho(\| \hat{u} \|) + \Delta(\hat{x}) + \begin{bmatrix}
W w - \hat{W} \hat{w} \\
h_2(x) - H \hat{h}_2 (\hat{x})
\end{bmatrix}^T
\begin{bmatrix}
X^{11} & X^{12} \\
X^{21} & X^{22}
\end{bmatrix}
\begin{bmatrix}
W w - \hat{W} \hat{w} \\
h_2(x) - H \hat{h}_2 (\hat{x})
\end{bmatrix}. \nonumber %\label{storageFn}
\end{align}
\end{definition}
\noindent Here, we relaxed the definition of storage functions given in \cite{ZA3} to {\it practical} storage functions by allowing the upper bound on their derivative to include a nonnegative function $\Delta(\hat{x})$. The term $v(x, \hat{x}, \hat{u})$ acts as the associated interface in Definition \ref{storageFnDef} by providing the concrete control input $u$. \added{We note that the purpose of matrix $H$ is to allow comparison between $h_2(x)$ and $\hat{h}_2(\hat{x})$, which can have different output dimensions. Similarly, matrices $W$ and $\hat{W}$ allow comparison between $w$ and $\hat{w}$. The choice of matrices $X^{11}, X^{12}, X^{21}$, and $X^{22}$ specify the type of dissipativity property being described \cite{arcak2016networks}.}
% This contrasts with Definition 3.1 in \cite{zamani2017compositional}, where the input $u$ still depends on $x$, $\hat{x}$, and $\hat{u}$, but not via an explicit function $v$. 

Next, we define the class of interconnected control systems that we consider in this paper:
%- this is the same class considered in \cite{zamani2017compositional}.

\begin{definition} \label{interconnectionDefinition}
Consider $N$ control subsystems $\Sigma_i = (\R^{n_i}, \R^{m_i}, \R^{p_i}, f_i, \R^{q_{1i}}, \R^{q_{2i}}, h_{1i}, h_{2i})$, $i = 1, \dots, N$, and a static matrix $M$ of appropriate dimension describing the coupling of these subsystems. The \textit{interconnected control system} $\Sigma = (\R^n, \R^m, f, \R^q, h)$, denoted as $\mathcal{I}(\Sigma_1, \dots, \Sigma_N)$, is given by $n = \sum_{i = 1}^N n_i$, $m = \sum_{i=1}^N m_i$, $q = \sum_{i=1}^N q_{1i}$, and
\begin{align*}
&f(x, u) := [f_1(x_1, u_1, w_1); \dots; f_N(x_N, u_N, w_N)], \\
&h(x) := [h_{11}(x_1); \dots; h_{1N}(x_N)],
\end{align*}
where $u = [u_1; \dots; u_N] \in \mathbb{R}^n, \ x = [x_1; \dots; x_N] \in \mathbb{R}^m$, and with the internal variables constrained by
\begin{equation}\label{interconnection}
[w_1; \dots; w_N] = M[h_{21}(x_1); \dots; h_{2N}(x_N)].
\end{equation}
A depiction of an interconnected control system $\mathcal{I}(\Sigma_1, \dots, \Sigma_N)$ is given in Figure \ref{schematic}.
\begin{figure}[H]
\centering
\includegraphics[width = 0.325\columnwidth]{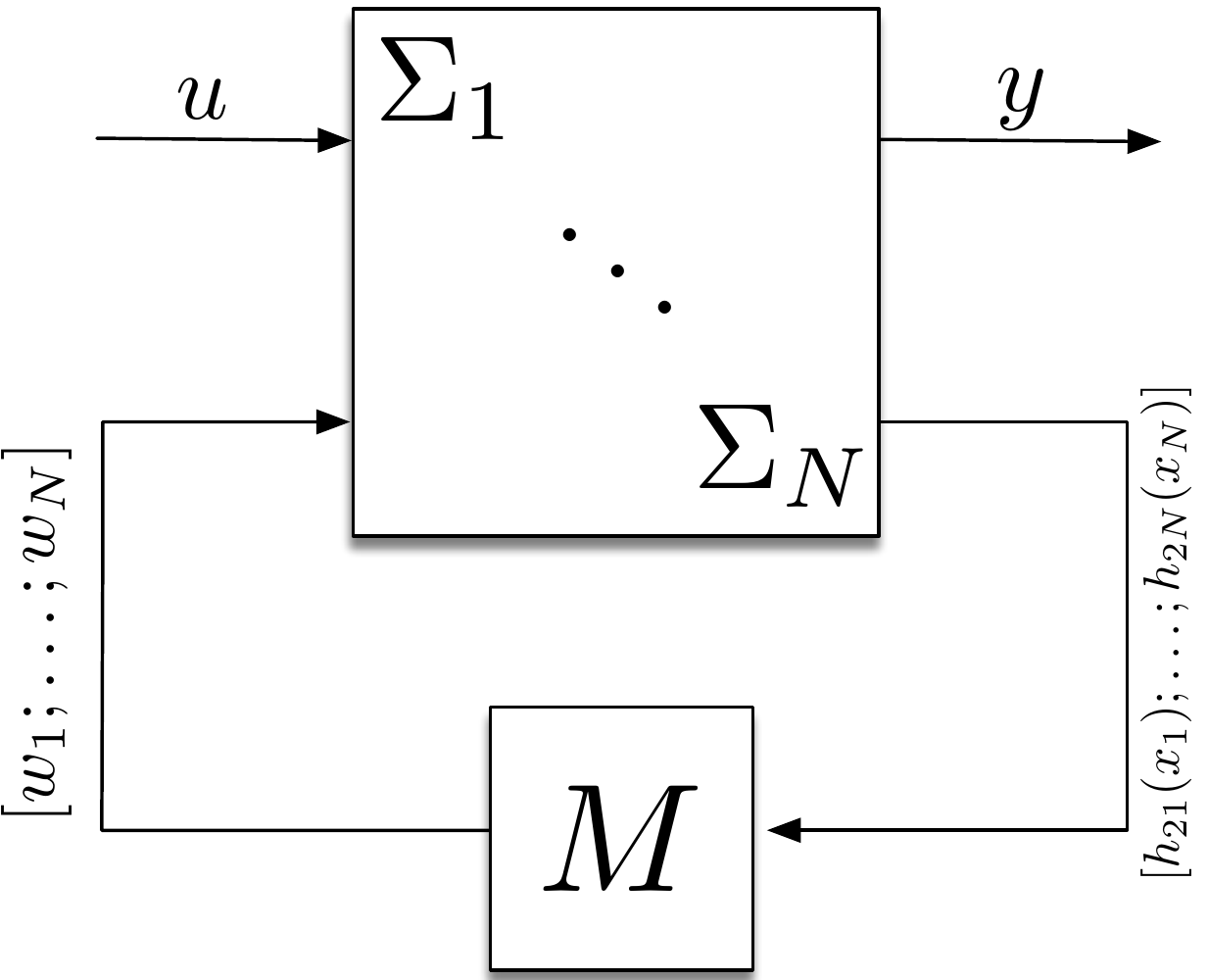}
\caption{An interconnection of $N$ control subsystems $\Sigma_1, \dots, \Sigma_N$.}
\label{schematic}
\end{figure}
\end{definition}

\subsection{Compositionality result}
We now provide a theorem containing our main result on the compositional construction of an abstraction and corresponding practical simulation function. In Definition \ref{storageFnDef} we included a nonnegative term $\Delta(\hat{x})$, allowing one to construct abstractions at the subsystem level by utilizing a relaxation similar to what was done in Section \ref{monolithic}. Our next result is to show that a similar relaxation can also be made at the level of the interconnected control system. We first review a theorem from \cite{ZA3} that constructs simulation functions from storage functions associated to subsystems; we then present a modified version with relaxed conditions.

\begin{theorem}\cite[Theorem 4.2]{ZA3} \label{compositionalityResult}
Consider the interconnected control system $\mathcal{I}(\Sigma_1, \dots, \Sigma_N)$ induced by $N$ control subsystems $\Sigma_i$ and the coupling matrix $M$. Suppose each subsystem $\Sigma_i$ admits an abstraction $\hat{\Sigma}_i$ and corresponding storage function $V_i$, each with the associated functions and matrices $\nu_i$, $\eta_i$, $\rho_i$, $v_i$, $H_i$, $W_i$, $\hat{W}_i$, $X_i$, $X_i^{11}$, $X_i^{12}$, $X_i^{21}$, and $X_i^{22}$ appearing in Definition \ref{storageFnDef} (by dropping term $\Delta(\hat{x})$). If there exist scalars $\mu_i > 0$, $i = 1, \dots, N$, and matrix $\hat{M}$ of appropriate dimension such that the following matrix (in)equality constraints
\begin{align}
 \begin{bmatrix}
W M \\
I_{\tilde{q}}
\end{bmatrix}^T
X( \mu_1 X_1, \dots, \mu_N X_N)
\begin{bmatrix}
W M \\
I_{\tilde{q}}
\end{bmatrix} &\leq 0, \label{equitableLMI} \\
 W M H &= \hat{W} \hat{M},
\label{equitabilityconstr}
\end{align}
are satisfied, where $\tilde{q} = \sum_{i = 1}^N q_{2i}$ and
\begin{small}
\begin{align}
& W := \text{diag}(W_1, \dots, W_N), \ \hat{W} := \text{diag}(\hat{W}_1, \dots, \hat{W}_N), \ H := \text{diag}(H_1, \dots, H_N), \label{blkMatrices} \\
& X(\mu_1 X_1, \dots, \mu_N X_N) := \begin{bmatrix}
\mu_1 X_1^{11} & & & \mu_1 X_1^{12} & & \\
& \ddots & & & \ddots & \\
& & \mu_N X_N^{11} & & & \mu_N X_N^{12} \\
\mu_1 X_1^{21} & & & \mu_1 X_1^{22} & & \\
& \ddots & & & \ddots & \\
& & \mu_N X_N^{21} & & & \mu_N X_N^{22}
\end{bmatrix}, \label{Xdissipativity}
\end{align}
\end{small} then
\begin{equation} \label{compositionalSimFun}
V(x, \hat{x}) := \sum_{i=1}^N \mu_i V_i(x_i, \hat{x}_i)
\end{equation}
is a simulation function from the interconnected control system $\hat{\Sigma} = \mathcal{I}(\hat{\Sigma}_1, \dots, \hat{\Sigma}_N)$, with the coupling matrix $\hat{M}$, to $\Sigma$. 
\end{theorem}
The following theorem relaxes (\ref{equitabilityconstr}) in Theorem \ref{compositionalityResult} as follows:
\begin{theorem}\label{compositionalityResultmod}
Suppose, instead of (\ref{equitabilityconstr}), one can only find a matrix $\hat{M}$ yielding a residual
\begin{equation} \label{compositionalResidual}
Y := W M H - \hat{W} \hat{M}
\end{equation}
which is nonzero, and all other hypotheses of Theorem \ref{compositionalityResult} hold with each $V_i$ being a practical storage function as in Definition \ref{storageFnDef}. Then \eqref{compositionalSimFun} is a practical simulation function from $\hat{\Sigma}$ to $\Sigma$ if there exist $\mu_i > 0$ and matrix $Z = Z^T \geq 0$ of appropriate dimensions such that the following matrix inequality constraint holds
%\begin{align} \label{relaxedCompositionalLMI}
%\begin{bmatrix}
%Y^T X^{11} Y - Z & Y^T X^{12} + Y^T X^{11} W M \\
%X^{21} Y + M^T W^T X^{11} Y & \begin{bmatrix}
%W M \\
%I_{\tilde{q}}
%\end{bmatrix}^T
%X
%\begin{bmatrix}
%W M \\
%I_{\tilde{q}}
%\end{bmatrix}
%\end{bmatrix} \nonumber \\[7.5pt]
%:= Q(Z, \mu_1, \dots, \mu_N) \leq 0,
%\end{align} where, to improve readability, we have written \eqref{Xdissipativity} as $X$, $\text{diag}(\mu_1 X_1^{11}, \dots, \mu_N X_N^{11})$ as $X^{11}$, and so on. In particular, the function $\Delta(\hat{x})$ in Definition \ref{simFunDef} is given by
\begin{equation}
Q(Z, \mu_1, \dots, \mu_N):=
\begin{bmatrix}
Y & W M \\
0 & I_{\tilde{q}}
\end{bmatrix}^T
X(\mu_1 X_1, \dots, \mu_N X_N)
\begin{bmatrix}
Y & W M \\
0 & I_{\tilde{q}}
\end{bmatrix} - \begin{bmatrix}
Z & 0 \\
0 & 0
\end{bmatrix} \le 0. \label{relaxedCompositionalLMI}
\end{equation}
In particular, the function $\Delta(\hat{x})$ in Definition \ref{simFunDef} is given by
\begin{equation} \label{nonequitableDeviation}
\Delta(\hat{x}) :=
\begin{bmatrix}
\hat{h}_{21}(\hat{x}_1) \\
\vdots \\
\hat{h}_{2N}(\hat{x}_N)
\end{bmatrix}^T Z
\begin{bmatrix}
\hat{h}_{21}(\hat{x}_1) \\
\vdots \\
\hat{h}_{2N}(\hat{x}_N)
\end{bmatrix} + \sum_{i = 1}^N \mu_i \Delta_i(\hat{x}_i).
\end{equation}
\end{theorem}
%\Majid{Very long theorem! We should break it in two theorems which the first one is just the recall from \cite{zamani2017compositional}.}

Theorem \ref{compositionalityResultmod} dropped the constraint \eqref{equitabilityconstr} from Theorem \ref{compositionalityResult}, resulting in a residual term \eqref{compositionalResidual}. The effect of this relaxation is then quantified via the term $\Delta(\hat{x})$, which is parameterized by the matrix $Z$ and scalars $\mu_i$ in \eqref{nonequitableDeviation}. Therefore, Theorem \ref{compositionalityResultmod} is beneficial when no matrix $\hat{M}$ satisfying \eqref{equitabilityconstr} exists. For such a scenario, we provide two optimization problems that can be solved in sequence to minimize the resulting $\Delta(\hat{x})$. First, with matrices $W$, $M$, $H$, and $\hat{W}$ fixed, we select the matrix $\hat{M}$ to minimize the residual \eqref{compositionalResidual} as measured by the Frobenius norm:\\
%. This task is outlined in the following optimization problem. \\

\noindent \textbf{Optimization Problem 2:}
\begin{align}
\text{minimize } & \|W M H - \hat{W} \hat{M}\|_F \quad \text{over} \ \hat{M}. \nonumber
\end{align}

With $\hat{M}$ thus selected, our next goal is to find a minimal $\Delta(\hat{x})$ as defined in \eqref{nonequitableDeviation}. We first introduce a diagonal scaling matrix $S$ that induces the functions $\tilde{h}_{2i}$, $i = 1, \dots, N$, as follows
\begin{equation}
\begin{bmatrix}
\tilde{h}_{21}(\hat{x}_1) \\
\vdots \\
\tilde{h}_{2N}(\hat{x}_N)
\end{bmatrix} := 
\underbrace{\begin{bmatrix}
s_1 I_{\hat{q}_{21}} & & \\
& \ddots & \\
& & s_N I_{\hat{q}_{2N}} 
\end{bmatrix}}_{:= S}
\begin{bmatrix}
\hat{h}_{21}(\hat{x}_1) \\
\vdots \\
\hat{h}_{2N}(\hat{x}_N)
\end{bmatrix}. \nonumber
\end{equation}
In particular, the scalars $s_i > 0$ are to be chosen so the outputs of the functions $\tilde{h}_{2i}(\hat{x}_i)$, $i=1,\dots,N$, \added{are comparable in order of magnitude}. \addedTwo{Next, we define the scalars $r_i > 0$, $i = 1, \dots, N$, which scale the functions $\Delta_i(\hat{x}_i)$ in the same way}. Then, we propose finding a minimal $\Delta(\hat{x})$ by solving the following optimization problem. \\

\noindent \textbf{Optimization Problem 3:}
\begin{align}
\text{minimize }& \quad \textbf{tr} (S^{-T} Z S^{-1}) + \sum_{i = 1}^N \frac{\mu_i}{r_i} \quad \text{over } Z \geq 0,  \ \mu_i \geq 1, \ i = 1, \dots, N, \nonumber \\[5pt]
\text {subject to }& \quad Q(Z, \mu_1, \dots, \mu_N) \leq 0. \label{lmiconsr}
\end{align}
\addedTwo{In particular, here the objective function represents our goal of minimizing $\Delta(\hat{x})$ in \eqref{nonequitableDeviation}, thus minimizing the error bound obtained via Theorem 1}. Here, we constrain $\mu_i \geq 1$ so that the \added{decision variables $\mu_i$ and $Z$ do not become too small and, as a result, poorly scaled}. We note that Optimization Problems 2 and 3 are both conic, and thus can be solved with a conic optimization tool such as MOSEK \cite{mosek}.

% \Majid{Should not we add small details on how to solve those optimization problems? E.g. via semidefinite programming ...}

\section{Aggregation} \label{aggregationSection}

A common approach to model order reduction in large scale systems is {aggregation}, which combines physical variables into a small number of groups and studies the interaction among these groups.  Examples include power systems, where geographical areas in which generators swing in synchrony are aggregated into {\it equivalent machines} \cite{chow}, and multicellular ensembles, where groups of cells exhibiting homogeneous behavior are represented with lumped biochemical reaction models \cite{FerArc13}.  

In this section we study a network of agents and first review an {\it equitable partition} criterion for aggregation when the agents have identical models.  We next relax the identical model assumption and the equitability criterion by using the results of the previous sections.  We formulate an optimization problem that penalizes the violation of the equitability condition when partitioning the agents into aggregate groups and, finally, study a special class of systems that encompasses the temperature control example in the next section.

\subsection{Equitable partition criterion for aggregation}
Consider $L$ agents with identical dynamical models:
\begin{align}\label{agent}
\dot{\xi}^\ell (t) &=g(\xi^\ell(t),\upsilon^\ell(t),\omega^\ell(t))\\
\zeta_1^\ell(t)&=\varsigma(\xi^\ell(t)) \label{readout}
\\
\label{agentout}
\zeta_2^\ell(t) &= \sigma(\xi^\ell(t)) \qquad \ell=1,2,\cdots,L,
\end{align}
 $\xi^\ell(t) \in \mathbb{R}^n$, $\upsilon^\ell(t) \in \mathbb{R}^{m}$, $\omega^\ell(t) \in \mathbb{R}^p$, $\zeta_1^\ell(t) \in \mathbb{R}^q$, $\zeta_2^\ell(t) \in \mathbb{R}^p$, for any $t\geq0$, interconnected according to the relation
\begin{equation} \label{intercon}
\begin{bmatrix} \omega^1(t) \\ \vdots \\ \omega^L(t) \end{bmatrix} =(\tilde{M}\otimes I_p) \begin{bmatrix} \zeta_2^1(t) \\ \vdots \\ \zeta_2^L(t) \end{bmatrix}, \quad \tilde{M}\in \mathbb{R}^{L\times L}.
\end{equation}
We partition the agents $\{1,\dots,L\}$ into $N\le L$ groups and describe the assignment of the agents to the groups with the $L\times N$ partition matrix 
\begin{equation}\label{Pdef1}
P_{\ell,i}=\left\{ \begin{array}{ll}  1 & \mbox{if} \ \ell \, \in \, \mbox{\small group}\ i \\ 0 & \mbox{otherwise.} \end{array}\right.
\end{equation}
\added{In particular, each agent is assigned to exactly one group, and each group must have at least one agent assigned to it}. We then aggregate the agents comprising each group into a single agent model that describes homogeneous behavior within the group. Thus, the abstraction for group $i$ is
\begin{align}\label{abs}
\dot{\hat{\xi}}_i(t)&=g(\hat{\xi}_i(t),\hat{\upsilon}_i(t),\hat{\omega}_i(t)) \\
\hat{\zeta}_{1i}(t)&=\mathbf{1}_{L_i}\otimes\varsigma(\hat{\xi}_i(t))
\\
\label{absout}
\hat{\zeta}_{2i}(t)&= \sigma(\hat{\xi}_i(t)) \qquad i=1,2,\cdots,N,
\end{align}
where $L_i$ is the number of agents in group $i$, $\hat{\xi}_i(t) \in \mathbb{R}^n$, $\hat{\upsilon}_i(t) \in \mathbb{R}^{m}$, $\hat{\omega}_i(t) \in \mathbb{R}^p$, $\hat{\zeta}_{1i}(t) \in \mathbb{R}^{qL_i}$, $\hat{\zeta}_{2i}(t) \in \mathbb{R}^p$, for any $t\geq0$, and the interconnection relation is
\begin{equation}\label{abs_intercon}
\begin{bmatrix} \hat{\omega}_1(t) \\ \vdots \\ \hat{\omega}_N(t) \end{bmatrix} =(\bar{M}\otimes I_p) \begin{bmatrix} \hat{\zeta}_{21}(t) \\ \vdots \\ \hat{\zeta}_{2N}(t) \end{bmatrix},
\end{equation}
where  $\bar{M}\in \mathbb{R}^{N\times N}$ is to be selected.

For the groups to exhibit perfectly homogeneous behavior, the trajectories must converge to and remain on the subspace where $\xi^\ell=\hat{\xi}_i$ for each $\ell$ in group $i$, $i=1,\dots,N$.  The invariance of this subspace is ensured if  $\upsilon^\ell=\hat{\upsilon}_i$ and $\omega^\ell=\hat{\omega}_i$ on the subspace, because $\xi^\ell(0)=\hat{\xi}_i(0)$, $\upsilon^\ell=\hat{\upsilon}_i$ and $\omega^\ell=\hat{\omega}_i$ imply $\dot{\xi}^\ell=\dot{\hat{\xi}}_i$
 by (\ref{agent}) and (\ref{abs}).  The internal inputs $\omega^\ell$, however, are not independent variables and the condition that $\omega^\ell=\hat{\omega}_i$ for having $\xi^\ell=\hat{\xi}_i$ for each $\ell$ in group $i$ must be further examined. To do so, first note from (\ref{agentout}) and (\ref{absout}) that $\xi^\ell=\hat{\xi}_i$ implies $\zeta_2^\ell=\hat{\zeta}_{2i}$,  which means 
\begin{equation}\label{identical_out}
\begin{bmatrix} \zeta_2^1(t) \\ \vdots \\ \zeta_2^L(t) \end{bmatrix} = (P\otimes I_p) \begin{bmatrix} \hat{\zeta}_{21}(t) \\ \vdots \\ \hat{\zeta}_{2N}(t) \end{bmatrix} \nonumber
\end{equation}
and, from (\ref{intercon}),
  \begin{equation}\label{w_from_haty}
\begin{bmatrix} \omega^1(t) \\ \vdots \\ \omega^L(t) \end{bmatrix} =(\tilde{M}P\otimes I_p) \begin{bmatrix} \hat{\zeta}_{21}(t) \\ \vdots \\ \hat{\zeta}_{2N}(t) \end{bmatrix}.
\end{equation}
The desired condition is  $\omega^\ell(t)=\hat{\omega}_i(t)$  for each $\ell$ in group $i$, that is
  \begin{equation}\label{identical_in}
\begin{bmatrix} \omega^1(t) \\ \vdots \\ \omega^L(t) \end{bmatrix} = (P\otimes I_p) \begin{bmatrix} \hat{\omega}_1(t) \\ \vdots \\ \hat{\omega}_N(t) \end{bmatrix}, \nonumber
\end{equation}
which is consistent with (\ref{w_from_haty}) if and only if $\bar{M}$ in (\ref{abs_intercon}) satisfies
 \begin{equation} \label{equit}
 \tilde{M}P=P\bar{M}.
 \end{equation}
  Thus, the invariance of the subspace $\xi^\ell=\hat{\xi}_i$ for each $\ell$ in group $i$ hinges upon the property (\ref{equit}), formalized in the following definition:
  
 \begin{definition} 
 Given $L$ agents with interconnection matrix $\tilde{M}\in \mathbb{R}^{L \times L}$, a partition into $N$ groups is said to be {\it equitable} if the partition matrix $P$ in (\ref{Pdef1}) satisfies (\ref{equit}) 
for some $\bar{M}\in \mathbb{R}^{N\times N}$.
\end{definition}

\begin{figure}
\centering
\includegraphics[width = 0.45\columnwidth]{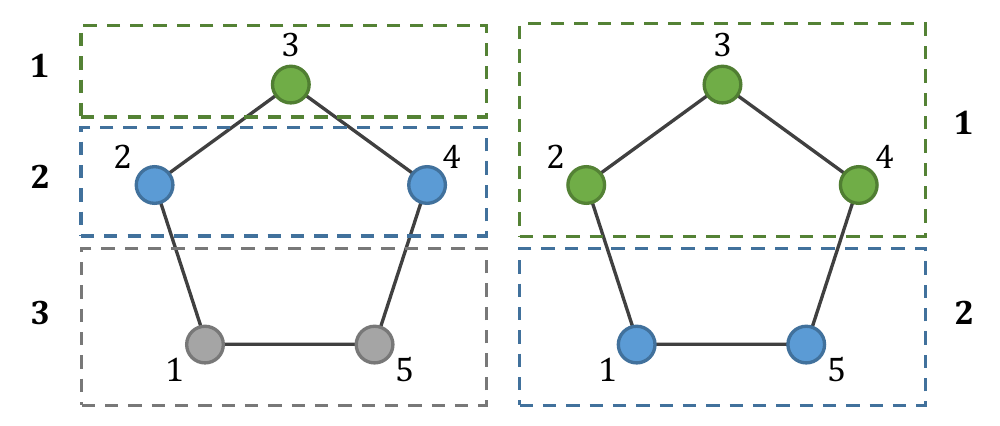}
%\includegraphics[width = 0.45\columnwidth]{images/rooms_graph.pdf}
% \vspace{-5mm}
\caption{\added{An equitable partition of a circle graph with $L=5$ nodes into three groups (left). Note that the partition into two groups (right) is \emph{not} equitable.}
\label{roomsGraph}}
\end{figure}

To provide intuition behind equitability, suppose $\tilde{M}$ corresponds to the Laplacian matrix of an unweighted, undirected graph, where each node represents an agent and edges are drawn between agents which are connected to one another. \added{In this case, a partition of the graph is equitable if each node in group $k$ has exactly $c_{k \ell}$ neighbors in group $\ell$, regardless of which node in class $k$ we select \cite{godsil2013algebraic}. Here, the constant $c_{k \ell}$ depends on $k$ and $\ell$.} As an illustration, an equitable partition of a five-node circle graph is displayed in Figure \ref{roomsGraph} (left), where group $1$ consists of node $3$, group $2$ consists of nodes $\{2,4\}$, and group 3 consists of nodes $\{1,5\}$. Each node in group $2$ is connected to \added{$c_{21} = 1$} node in group $1$ and \added{$c_{23} = 1$} node in group $3$. On the other hand, note that the partition displayed in Figure \ref{roomsGraph} (right) into groups $\{2, 3, 4\}$ and $\{1, 5\}$ is \textit{not} equitable. \added{Although we discussed unweighted graphs for simplicity, the extension to weighted graphs is straightforward by considering the sum of the edge weights connected to a particular node instead of the number of neighbors.}

\subsection{Relaxing the identical agent and equitable partition assumptions}
The assumptions that the agent dynamics be identical and that an equitable partition exist for their interconnection can be restrictive in practice.  The control specifications may further limit the choice of partition, since the states of agents in the same group are lumped together in the abstraction and the specifications cannot distinguish between them.  

Here we relax both assumptions using the results of Section \ref{Composition-Section}.  First we replace the agent dynamics (\ref{agent}) with
\begin{equation}\label{agent2}
\dot{\xi}^\ell(t) = g^\ell(\xi^\ell(t),\upsilon^\ell(t),\omega^\ell(t)),
\end{equation}
where $g^\ell: \mathbb{R}^n\times \mathbb{R}^m \times \mathbb{R}^p\rightarrow \mathbb{R}^n$, $\ell=1,\dots,L$, allow for deviations from the nominal model $g$ used in the abstraction (\ref{abs}). \added{We note that it is also possible to define separate nominal dynamics for each group in order to minimize said deviations further. For simplicity, here we use the same nominal model for each group.} In preparation for constructing a simulation function, we assume that there exist practical storage functions from the agents to the nominal model with identical supply rates as the following:

\begin{assumption} \label{A1} There exist a locally Lipschitz function $\tilde{v}^\ell: \mathbb{R}^n\times \mathbb{R}^n \times \mathbb{R}^m \rightarrow \mathbb{R}^m$, a continuously differentiable function $\tilde{V}^\ell: \mathbb{R}^n\times \mathbb{R}^n  \rightarrow \mathbb{R}_{\ge 0}$, $\tilde{\nu}^\ell,\tilde{\eta}^\ell\in \mathcal{K}_\infty$, $\tilde{\rho}^\ell\in \mathcal{K} \cup \{0\}$,  $\tilde{\Delta}^\ell:\mathbb{R}^n\rightarrow \mathbb{R}_{\ge 0}$, and a matrix $\tilde{X}=\tilde{X}^T\in \mathbb{R}^{2p\times 2p}$ such that for all $x \in \mathbb{R}^n$, $\hat{x} \in \mathbb{R}^n$, $\hat{u} \in \mathbb{R}^{m}$, $w \in \mathbb{R}^p$, $\hat{w} \in \mathbb{R}^p$,
\begin{equation}  \label{lowerbound}
 \tilde{\nu}^\ell(\| \varsigma(x)-\varsigma(\hat{x}) \|) \le \tilde{V}^\ell(x,\hat{x}),
 \end{equation}
\begin{align}\label{derivative}
\frac{\partial \tilde{V}^\ell(x,\hat{x})}{\partial x}g^\ell&(x,\tilde{v}^\ell(x,\hat{x},\hat{u}),{w})+\frac{\partial \tilde{V}^\ell(x,\hat{x})}{\partial \hat{x}}g(\hat{x},\hat{u},\hat{w})\\
 \le& -\tilde{\eta}^\ell(\tilde{V}^\ell(x,\hat x))+\tilde{\rho}^\ell(\|\hat{u}\|)+\tilde\Delta^\ell(\hat{x}) +\begin{bmatrix} w-\hat{w} \\ \sigma(x)-\sigma(\hat{x}) \end{bmatrix}^T \!\!\! \tilde{X} \! \begin{bmatrix} w-\hat{w} \\ \sigma(x)-\sigma(\hat{x}) \end{bmatrix}\!. \label{agent_dissip}
\end{align}
\end{assumption}
%\Majid{Shouldn't \eqref{lowerbound} be based on outputs of states?}

 In the next subsection we show a class of systems that satisfy Assumption \ref{A1}. One will see, in particular, that the term $\tilde\Delta^\ell(\hat{x})$ in (\ref{agent_dissip}) is critical for absorbing the mismatch between $g^\ell$ and $g$, which is due to the heterogeneity of the agent models. \added{In the example in Section 6, we show how the interface function can be used to help to shrink $\tilde\Delta^\ell(\hat{x})$ and satisfy Assumption 1 with a tight upper bound in \eqref{derivative}.}

We let each group $i=1,\dots,N$ in the partition define a subsystem,  and derive a composite storage function and dissipation inequality from Assumption \ref{A1}. 
Let $L_i\ge 1$ denote the number of agents in group $i$, $L_1+\cdots+L_N=L$, and define the state vector $x_i \in \mathbb{R}^{L_in}$ by concatenating the state vectors $x^\ell$ of the agents assigned to group $i$.
  Defining $u_i\in \mathbb{R}^{L_im}$, $w_i\in \mathbb{R}^{L_ip}$, $y_{2i}\in \mathbb{R}^{L_iq}$  and  $y_{2i}\in \mathbb{R}^{L_ip}$, we write the model for subsystem $i$ as
\begin{align}\label{subsys}
\dot{\xi}_i(t)&=f_i(\xi_i(t),\upsilon_i(t),\omega_i(t)) \\
\zeta_{1i}(t)&= h_{1i}(\xi_i(t))
\\
\label{subsysout}
\zeta_{2i}(t)&= h_{2i}(\xi_i(t)) 
\end{align}
where $f_i(\xi_i(t),v_i(t),\omega_i(t))$, $h_{1i}(\xi_i(t))$ and $h_{2i}(\xi_i(t))$ are obtained by concatenating $g^\ell(\xi^\ell(t),v^\ell(t),\omega^\ell(t))$, $\varsigma(\xi^\ell(t))$ and $\sigma(\xi^\ell(t))$, respectively, over each $\ell$ in group $i$. 

We assume, without loss of generality, that the agents are indexed such that the first $L_1$ constitute group $1$, the next $L_2$ group $2$, and so on. It then follows from (\ref{intercon}) that
\begin{equation}\label{intercon2}
\begin{bmatrix} \omega_1(t) \\ \vdots \\ \omega_N(t) \end{bmatrix} =(\tilde{M}\otimes I_p) \begin{bmatrix} \zeta_{21}(t) \\ \vdots \\ \zeta_{2N}(t) \end{bmatrix},
\end{equation}
since the respective vectors in (\ref{intercon}) and (\ref{intercon2}) are identical.  Without this assumption an appropriate permutation can be applied to the matrix $\tilde{M}$ and the subsequent results do not change.

Using Assumption \ref{A1} we let each agent $\ell$ in group $i$ apply the feedback $u^\ell=\tilde{v}^\ell(x^\ell,\hat{x}_i,\hat{u}_i)$, and define the practical storage function for subsystem $i$ to be
\begin{equation}\label{storage}
V_i(x_i,\hat{x}_i)=\sum_{\ell \, \in \, \mbox{\small group}\ i}\tilde{V}^\ell(x^{\ell},\hat{x}_i).
\end{equation}
Then, we obtain the dissipativity property:
\begin{align*}
\frac{\partial {V}_i(x_i,\hat{x}_i)}{\partial x_i}f_i&(x_i,u_i,w_i)\!+\!\frac{\partial {V}_i(x_i,\hat{x}_i)}{\partial \hat{x}_i}g(\hat{x}_i,\hat{u}_i,\hat{w}_i) \\
& = \sum_{\ell \, \in \, \mbox{\small group}\ i}\!\! \left\{\frac{\partial \tilde{V}^\ell(x^\ell,\hat{x}_i)}{\partial x^\ell}g^\ell(x^\ell,\tilde{v}^\ell(x^\ell,\hat{x}_i,\hat{u}_i),{w}^\ell)+\frac{\partial \tilde{V}^\ell(x^\ell,\hat{x}_i)}{\partial \hat{x}_i}g(\hat{x}_i,\hat{u}_i,\hat{w}_i)\right\} \\
&\le \sum_{\ell \, \in \, \mbox{\small group}\ i} \!\! \left\{-\tilde\eta^\ell(\tilde{V}^\ell(x^{\ell},\hat{x}_i)) + \tilde{\rho}^\ell(\|\hat{u}_i\|)+\tilde{\Delta}^\ell(\hat{x}_i)+\begin{bmatrix} w^\ell-\hat{w}_i \\ \sigma(x^\ell)-\sigma(\hat{x}_i) \end{bmatrix}^T \!\!\! \tilde{X} \! \begin{bmatrix} w^\ell-\hat{w}_i \\ \sigma(x^\ell)-\sigma(\hat{x}_i) \end{bmatrix}\right\} \\
& \le -\eta_{i}(V_i(x_i,\hat{x}_i))+\rho_{i}(\|\hat{u}_i\|)+\Delta_i(\hat{x}_i) +\begin{bmatrix} w_i-(\mathbf{1}_{L_i}\otimes I_p)\hat{w}_i \\ y_{2i}-(\mathbf{1}_{L_i}\otimes I_p)\hat{y}_{2i}\end{bmatrix}^T \!\!\!\!\! {X}_i \!\! \begin{bmatrix} w_i-(\mathbf{1}_{L_i}\otimes I_p)\hat{w}_i \\ y_{2i}-(\mathbf{1}_{L_i}\otimes I_p)\hat{y}_{2i}\end{bmatrix}
\end{align*} 
where, \added{for $s \in \R_{\geq 0}$ and $y \in \R^n$}, we define
\begin{align}\label{rhodelta}
\eta_{i}(s) \! &:= \!\! \min_{z\in \mathbb{R}^{L}_{\ge 0}} \! \sum_{\ell \, \in \, \mbox{\small group}\ i}\!\!\!\!\!\! \tilde{\eta}^\ell(z_\ell)   \ \mbox{s.t.} \!\!\!\!\!\!\! \sum_{\ell \, \in \, \mbox{\small group}\ i}\!\!\!\!\!\!\!\!z_\ell=s, \quad \rho_i(s) \!:= \!\!\!\!\! \sum_{\ell \, \in \, \mbox{\small group}\ i}\!\!\!\!\!\!\!  \tilde{\rho}^\ell(s), \quad \Delta_i(y)\! := \!\!\!\!\!  \sum_{\ell \, \in \, \mbox{\small group}\ i}\!\!\!\!\!\!\!  \tilde{\Delta}^\ell(y),
\\
\label{Xi}
X_{i} &:=\begin{bmatrix}I_{L_i}\otimes \tilde{X}^{11} & I_{L_i}\otimes \tilde{X}^{12} \\ I_{L_i}\otimes \tilde{X}^{21} & I_{L_i}\otimes \tilde{X}^{22}\end{bmatrix},
\end{align}
and where $\tilde{X}^{11}$, $\tilde{X}^{12}$, $\tilde{X}^{21}$, $\tilde{X}^{22}$ denote $p\times p$ matrices obtained by partitioning $\tilde{X}\in \mathbb{R}^{2p\times 2p}$ conformally.
Defining, in addition,
\begin{equation}\label{WH}
W_i:=I_{L_ip}, \ \hat{W}_i=H_i:=\mathbf{1}_{L_i}\otimes I_p
\end{equation}
and
\begin{equation}\label{nu}
 \nu_{i}(s):=\min_{z\in \mathbb{R}^{L}_{\ge 0}}  \sum_{\ell \, \in \, \mbox{\small group}\ i}\!\!\!\! \tilde{\nu}^\ell(z_\ell)   \quad \mbox{s.t.} \sum_{\ell \, \in \, \mbox{\small group}\ i}\!\!\!\!  z_\ell=s,
 \end{equation}
 we summarize the conclusion in the following proposition:

 \begin{proposition} \label{Prop1} Suppose the agents $\ell=1,\dots,L$ satisfy Assumption \ref{A1}, and each subsystem $i=1,\dots,N$ is defined as in (\ref{subsys})-(\ref{subsysout}), with the abstraction (\ref{abs})-(\ref{absout}) obtained 
 by aggregating $L_i$ agents.
 Then $V_i$ in \eqref{storage} is a practical storage function
as in Definition \ref{storageFnDef}, with (\ref{rhodelta})-(\ref{nu}), $\hat{h}_{1i}(\hat{x}_i)=\mathbf{1}_{L_i}\otimes \varsigma(\hat{x}_i)$, and $\hat{h}_{2i}(\hat{x}_{i})=\sigma(\hat{x}_{i})$.
\end{proposition}

We next examine the conditions of Theorem \ref{compositionalityResult} and Theorem \ref{compositionalityResultmod}. From (\ref{WH}) and (\ref{blkMatrices}) we have:
\begin{equation}\label{Wdef}
 W= I_{Lp}, \ \mbox{and} \ 
 \hat{W}=H=P\otimes I_p,
 \end{equation}
where
  \begin{equation}\label{Pdef2}
  P=\text{diag}(\mathbf{1}_{L_1}, \cdots, \mathbf{1}_{L_N})
  %P=\begin{bmatrix} \mathbf{1}_{L_1} & & \\ & \ddots & \\ & & \mathbf{1}_{L_N}\end{bmatrix}.
\end{equation}
Since we assumed that the agents are indexed such that the first $L_1$ constitute group $1$, the next $L_2$ group $2$, and so on, the definition of $P$ in (\ref{Pdef2}) is consistent with the partition matrix defined in (\ref{Pdef1}).
If the subsystem abstractions are interconnected as in (\ref{abs_intercon}), then  $\hat{M}=\bar{M}\otimes I_p$ and, thus, condition (\ref{equitabilityconstr}) of Theorem \ref{compositionalityResult} is identical  to the equitability criterion (\ref{equit}). This means that we can relax the equitability condition with Theorem \ref{compositionalityResultmod}.  The first residual term  in (\ref{nonequitableDeviation}) is then due to the relaxation of equitability, and the second term is due to model variations of non-identical agents, absorbed into $\tilde\Delta^\ell$ in Assumption \ref{A1} and combined into 
 $\Delta_i$ in (\ref{rhodelta}).

\subsection{An optimization problem for near-equitability}
We note that relaxing the equitability condition \eqref{equit} results in a residual term given by
\begin{equation} \label{equitRelaxation}
\bar{Y} := \tilde{M} P - P \bar{M}.
\end{equation}
Our goal now becomes choosing a partition of the agents - equivalently, a partition matrix $P$ and coupling matrix $\bar{M}$ - such that \eqref{equitRelaxation} is minimized. We propose approaching this task in two steps. First, we allow for some agents to be assigned to groups by hand. \added{Since aggregated agents share the same specification, this allows one to assign agents to separate groups if they require separate specifications. Conversely, one can also assign agents to the same group if it is desirable for them to abide by the same specification}. In the second step, the remaining agents are to be assigned to groups automatically via an optimization problem to be defined next. The pre-assigned agents induce an $L \times N$ matrix $\bar{P}$ as follows
\begin{equation}\label{preassign1}
\bar{P}_{\ell,i}=\left\{ \begin{array}{ll}  1 & \mbox{if $\ell$ is pre-assigned to group}\ i \\ 0 & \mbox{otherwise} \end{array}\right.
\end{equation}
as well as a diagonal matrix
\begin{equation} \label{preassign2}
T =\text{diag}(t_1, \cdots, t_N)
%\begin{bmatrix}
%t_1 & & \\
%& \ddots & \\
%& & t_N
%\end{bmatrix}
\end{equation}
where $t_i$ is the number of agents pre-assigned to group $i$. We note that if an agent $\ell$ is not pre-assigned to any group, then the corresponding row $\ell$ of $\bar{P}$ will contain only zeros.

To partition the remaining agents, we solve a mixed-integer program. We model $\bar{M}$ as a continuous decision variable and, noting \eqref{Pdef1}, model $P$ as a binary decision variable. The objective function of our problem is the
Frobenius norm
of the residual term $\bar{Y}$, the minimization of which yields an equitable partition when one exists, and a near-equitable partition otherwise. 

We also note it is possible to enforce \eqref{equitRelaxation} using \emph{linear} constraints. Since $\tilde{M}$ is fixed, the term $\tilde{M} P$ is linear - the problematic term is $P \bar{M}$, as it is the product of two decision variables. Linearity is achieved with a reformulation, implemented as the command ``binmodel" \cite{binmodel} in the toolbox YALMIP \cite{lofberg2004yalmip}. To see the idea for the scalar case, consider the product of a binary variable $p \in \{0, 1\}$ and a continuous variable $m \in \R$. Suppose that $m$ has lower bound $\underline{m} \in \R$ and upper bound $\overline{m} \in \R$.  Then, the product $p \cdot m$ can be replaced with a continuous auxiliary variable $y \in \R$ by including the following linear constraints
\begin{equation} \label{bigMconstraints}
\underline{m} p \leq y \leq \overline{m} p, \quad \underline{m} (1-p) \leq m - y \leq \overline{m} (1-p). \nonumber
\end{equation}
This procedure can be applied in a similar fashion to \eqref{equitRelaxation}. Thus, the following optimization problem can be cast as a mixed-integer quadratic program with linear constraints: \\

\noindent \textbf{Optimization Problem 4:} \label{opt2}
\begin{align}
\text{minimize } & \|\bar{Y}\|_F \quad \text{over} \ P, \ \bar{M} \nonumber \\[5pt]
\text {such that } & P \text{ is binary}, \label{eqconstr1} \\
& P \mathbf{1}_N = \mathbf{1}_L, \label{eqconstr2} \\
& \textbf{1}_L^T P \geq \mathbf{1}_N, \label{eqconstr3} \\
& \bar{P}^T P = T, \label{eqconstr4} \\
& \bar{Y} = \tilde{M} P - P \bar{M}, \label{eqconstr5}
\end{align}

\noindent where \eqref{eqconstr2} ensures each node is assigned to exactly one class, \eqref{eqconstr3} requires that each class has at least one node assigned to it, and \eqref{eqconstr4} assures that the pre-assignments represented by $\bar{P}$ and $T$, as defined in \eqref{preassign1} and \eqref{preassign2}, are respected. \added{We note that Optimization Problem 4 is a mixed integer quadratic program and therefore can be solved with an optimization tool such as Gurobi \cite{gurobi}.}

Note that Optimization Problem 4 minimizes the same residual as Optimization Problem 2, since $Y$ in \eqref{compositionalResidual} is equal to $\bar{Y} \otimes I_p$. However, here we have the additional flexibility of adjusting $P$, whereas the equivalent matrices $\hat{W}$ and $H$ in Optimization Problem 2 are fixed. Furthermore, since $\bar{M}$ is selected to minimize the Frobenius norm, the special structure of the matrix $P$ implies that $\bar{Y}$ has the following property:

\begin{lemma} \label{ybarLemma}
The matrix $\bar{Y}$ obtained by solving Optimization Problem 4 satisfies $\bar{Y}^T \mathbf{1}_L = 0$.
\end{lemma}

\noindent We will refer back to this fact after we state Theorem \ref{tempTheorem}, at which point it will become relevant. 

\subsection{A special class of agent models}
We now study a class of agent models of the form (\ref{readout}), (\ref{agentout}), (\ref{agent2}) with
\begin{equation}\label{specialclass}
g^\ell(x,u,w)=\alpha_\ell(x)+\beta_\ell(x)u+Bw, \ \varsigma(x)=x, \ \sigma(x)=Cx,
\end{equation}
where $\alpha_\ell: \mathbb{R}^n \rightarrow \mathbb{R}^n$ and $\beta_\ell: \mathbb{R}^n \rightarrow \mathbb{R}^{n\times m}$ are allowed to vary by agent $\ell$ and are replaced with nominal ones $\alpha: \mathbb{R}^n \rightarrow \mathbb{R}^n$ and $\beta: \mathbb{R}^n \rightarrow \mathbb{R}^{n\times m}$, respectively,
in the abstraction (\ref{abs})-(\ref{absout}):
\begin{equation}\label{specialabs}
g(\hat{x},\hat{u},\hat{w})=\alpha(\hat{x})+\beta(\hat{x})\hat{u}+B\hat{w}.
\end{equation}
\added{We note that $\alpha_\ell$ in \eqref{specialclass} is assumed to be continuously differentiable.} The following proposition gives sufficient conditions under which Assumption \ref{A1} holds for  (\ref{specialclass}) and (\ref{specialabs}) above:
\begin{proposition} \label{Prop2}
If there exists $\tilde{v}^\ell: \mathbb{R}^n\times \mathbb{R}^n \times \mathbb{R}^m \rightarrow \mathbb{R}^m$, $\tilde{\rho}^\ell\in \mathcal{K} \cup \{0\}$, constants $\lambda_\ell$, $\vartheta_\ell$, and $n\times n$ matrix $Q_\ell=Q^T_\ell>0$ such that, for all $x \in \mathbb{R}^n$, $\hat{x} \in \mathbb{R}^n$,
\begin{align} \label{Cond1}
& Q_\ell\left(\frac{\partial \alpha_\ell(x)}{\partial x}\right)\!+\!\left(\frac{\partial \alpha_\ell(x)}{\partial x}\right)^T\!\!\!\! Q_\ell\le 2\lambda_\ell I_n \\
 \label{Cond2}
& (x-\hat{x})^TQ_\ell\left(\beta_\ell(x)\tilde{v}^\ell(x,\hat{x},\hat{u})-\beta(\hat{x})\hat{u}\right)
\le \vartheta_\ell \|x-\hat{x}\|^2 +\tilde{\rho}^\ell(\|\hat{u}\|) \\
 \label{Cond3}
& \lambda_\ell+\vartheta_\ell<0 \\
 \label{Cond4}
& Q_\ell B=C^T,
\end{align}
then Assumption \ref{A1} holds with
\begin{align} \nonumber
\tilde{V}(x, \hat{x}) &= \frac{1}{2} (x - \hat{x})^T Q_\ell (x - \hat{x}), \quad \tilde{\eta}^\ell(s)=\frac{2\varepsilon_\ell}{\lambda_{\rm max}(Q_\ell)}s, \quad \tilde{\Delta}^\ell(\hat{x})=\frac{1}{4(|\lambda_\ell +\vartheta_\ell |-\varepsilon_\ell)}\|Q_\ell(\alpha_\ell(\hat{x})-\alpha(\hat{x}))\|^2, \\
\tilde{X}&=\frac{1}{2}\begin{bmatrix} 0 & I_p \\ I_p & 0\end{bmatrix} \label{passive}
\end{align}
for any choice of $\varepsilon_\ell \in (0, |\lambda_\ell+\vartheta_\ell|)$.
\end{proposition}
% \Majid{Can you elaborate more on conditions \eqref{Cond1}-\eqref{Cond4}? Which classes of systems satisfy those conditions (e.g. incrementally stabilizable or contractable?)}

Note that the conditions \eqref{Cond1} - \eqref{Cond4} imply that the system in \eqref{specialclass} is incrementally stabilizable.  We also note, in particular, that the term $\tilde{\Delta}^\ell(\hat{x})$ is due to the deviation of $\alpha_\ell(\hat{x})$ from $\alpha(\hat{x})$. Under the hypotheses of Proposition \ref{Prop2} it follows from Proposition \ref{Prop1} that the subsystems and their abstractions satisfy the dissipativity property in Definition \ref{storageFnDef} with 
$$
X_{i}=\frac{1}{2}\begin{bmatrix} 0 & I_{L_ip} \\ I_{L_ip} & 0\end{bmatrix}
$$
and, if we use identical weights $\mu_i=1$, $i=1,\dots,N,$ then the matrix $X$ in Theorem \ref{compositionalityResult} is
$$
{X}=\frac{1}{2}\begin{bmatrix} 0 & I_{Lp} \\ I_{Lp} & 0\end{bmatrix}.
$$
Since $W= I_{Lp}$ by (\ref{Wdef}),  condition (\ref{equitableLMI}) of Theorem \ref{compositionalityResult}
is
\begin{equation}
\begin{bmatrix} WM \\ I \end{bmatrix}^TX\begin{bmatrix} WM \\ I \end{bmatrix}
=\frac{1}{2}(\tilde{M}+\tilde{M}^T)\otimes I_p \le 0. \nonumber
\end{equation}

\begin{theorem} \label{tempTheorem} Suppose the agents $\ell=1,\dots,L$ are described by (\ref{readout}) - (\ref{intercon}), (\ref{agent2}),
with the special form (\ref{specialclass}) and interconnection matrix
\begin{equation}
\tilde{M}+\tilde{M}^T\le 0, \label{interconnectionCondition}
\end{equation}
and let the hypothesis of Proposition \ref{Prop2} hold. If the partition of the agents is equitable, then $V$ in \eqref{compositionalSimFun} is a practical simulation function from $\hat{\Sigma}$ to $\Sigma$ with $V_i$ as in \eqref{storage} and $\mu_i = 1$, $i = 1,\dots,N$. If the equitability condition \eqref{equit} is relaxed so $\bar{Y}$ in \eqref{equitRelaxation} is nonzero, then $V$ is a practical simulation function if there exists a matrix $Z = Z^T \geq 0$ satisfying \eqref{relaxedCompositionalLMI} with $Y = \bar{Y} \otimes I_p$ and $\mu_i = 1$ for $i = 1,\dots,N$. Furthermore, $\mathcal{N}(\tilde{M} + \tilde{M}^T) \subseteq \mathcal{N}(\bar{Y}^T)$ is a necessary and sufficient condition for such a $Z$ to exist.
\end{theorem}

The matrix $Z$ in Theorem \ref{tempTheorem} can be found by solving Optimization Problem 3, where we append the constraint $\mu_i = 1$ for $i = 1, \dots, N$. Furthermore, when $\bar{M}$ and $\bar{Y}$ are obtained via Optimization Problem 4, the null space condition of Theorem 5 holds automatically if $\mathcal{N}(\tilde{M} + \tilde{M}^T)$ is spanned by $\mathbf{1}_L$, since $\bar{Y}$ satisfies $\bar{Y}^T \mathbf{1}_L = 0$ from Lemma \ref{ybarLemma}. More generally, we also note if the stronger condition
\begin{equation}
\tilde{M} + \tilde{M}^T < 0 \nonumber
\end{equation}
on the interconnection matrix holds, then the null space condition is satisfied since $\mathcal{N}(\tilde{M} + \tilde{M}^T) = \{\mathbf{0}_L\}$.

\section{Example} \label{temperatureExample}
\subsection{Room Temperature Model}
We now consider a temperature control application adapted from \cite{girard2016safety}. Our goal is to control the temperature of $L$ rooms connected in a circle. We model the dynamics of the temperature $\xi_\ell(t) \in \R$ in room $\ell \in \{1, \dots, L\}$ as
\begin{align} \label{tempDynamics}
\dot{\xi}^\ell(t) &= a_\ell(T_e - \xi^\ell(t)) + b_\ell (T_h - \xi^{\ell}(t)) \upsilon^\ell(t) + \gamma \omega^\ell(t), \nonumber \\
\omega^\ell(t) &= \xi^{\ell+1}(t) + \xi^{\ell-1}(t) - 2 \xi^\ell(t),
\end{align}
where $a_\ell, b_\ell, \gamma \in \R_{>0}$ are conduction coefficients (where the former two may depend on room index), $T_e$ and $T_h$ are the temperatures of the external environment and room heater, respectively, and $\upsilon^\ell$ is a control input. Furthermore, we let $\xi^0 = \xi^L$ and $\xi^1 = \xi^{L+1}$ so that the indices in \eqref{tempDynamics} are valid for rooms $\ell = 1$ and $\ell = L$. Note that this model can be represented as in \eqref{specialclass} with $\alpha_\ell (s) = a_\ell (T_e - s)$, $\beta_\ell(s) = b_\ell (T_h - s)$, $B = \gamma$, and $C = 1$. Furthermore,  the coupling matrix is given by:
\begin{equation} \label{exampleLaplacian}
\tilde{M} = \begin{bmatrix}
-2 & 1 & 0 & \cdots & \cdots & 1 \\
1 & -2 & 1 & 0 & \cdots & 0 \\
0 & \ddots & \ddots & \ddots & \ddots & \vdots \\
\vdots & \ddots & \ddots & \ddots & \ddots & 0 \\
0 & & \ddots & \ddots & \ddots & 1\\
1 & 0 & \cdots& 0 & 1 & -2
\end{bmatrix}.
\end{equation}
%corresponding to the graph shown at the bottom of Figure \ref{roomsGraph}.

\begin{figure*}[tb]
\centering
\input{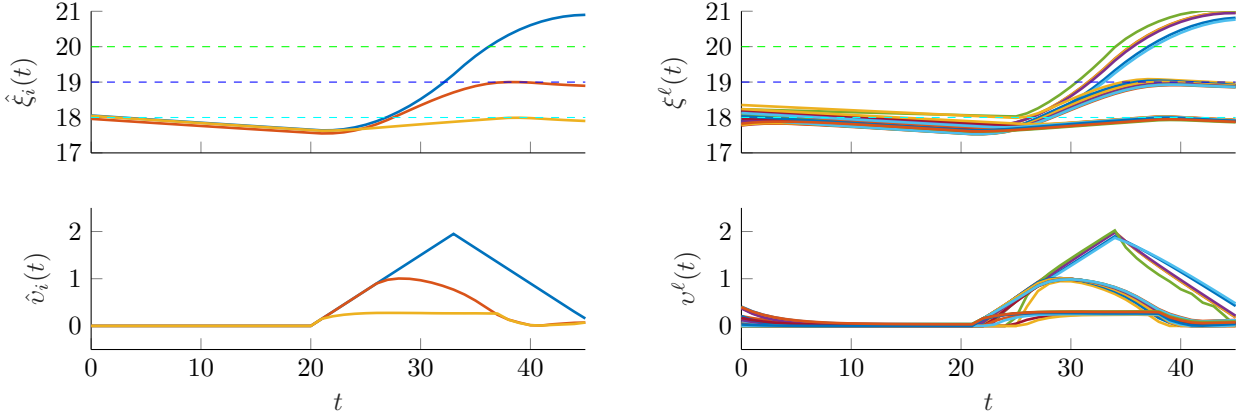}
\vspace{-2mm}
\caption{Simulation results for the temperature regulation example. We require the temperature in each area of the building to reach its corresponding target temperature range (indicated by the dashed lines) within 20 minutes after the signal is triggered. The signal is triggered at the 20 minute mark - the aggregate system (left) reaches the temperature target within 20 minutes, and the concrete system (right) closely follows the reference.} 
\label{temperatureSimulation}
\end{figure*}

\subsection{Aggregate Model}
For the aggregate model, we partition the rooms into $N \leq L$ distinct areas via Optimization Problem 4. The aggregate temperature $\hat{\xi}_i(t)$ in area $i \in \{ 1, \dots, N\}$ is governed by
\begin{equation}
\dot{\hat{\xi}}_i(t) = a (T_e - \hat{\xi}_i(t)) + b (T_h - \hat{\xi}_i(t)) \hat{\upsilon}_i(t) + \gamma \hat{\omega}_i(t) \nonumber
\end{equation}
where, in this case, the coupling $\hat{w}_i$ depends on the particular $\bar{M}$ we obtain by solving Optimization Problem 4. The conduction coefficients $a$ and $b$ in the nominal model are obtained by averaging over the conduction coefficients $a_\ell$ and $b_\ell$ for the individual rooms, so that $a := \frac{1}{L} \sum_{\ell = 1}^L a_\ell$ and $b := \frac{1}{L} \sum_{\ell = 1}^L b_\ell$. In this case, conditions \eqref{Cond1}, \eqref{Cond2}, and \eqref{Cond4} hold for the function
\begin{equation} \label{tempInterface}
\tilde{v}^\ell(x, \hat{x}, \hat{u}) = \frac{1}{b_\ell (T_h - x)} \left[ b (T_h - \hat{x}) \hat{u} - k_\ell (x - \hat{x}) \right]
\end{equation}
where $k_\ell \in \R_{\geq 0}$, $\tilde{\rho}^\ell( \| \hat{u} \|) = 0$, $\lambda_\ell = -a_\ell/\gamma$, $\vartheta_\ell = -k_\ell$, and $Q_\ell = 1/\gamma$. Furthermore, condition \eqref{Cond3} is satisfied if the gain $k_\ell$ is chosen such that $k_\ell > -a_{\ell} / \gamma$. Therefore, the result of Theorem \ref{tempTheorem} is applicable to this example, since $\tilde{M}= \tilde{M}^T \leq 0$. We also note that division by zero in \eqref{tempInterface} can be avoided by imposing constraints on $\hat{x}$ in a formal synthesis procedure - indeed, by combining this with a bound on the error between $x$ and $\hat{x}$, we can conclude that $x$ will never reach the heater temperature $T_h$. \added{A similar line of reasoning ensures the inputs of the aggregated systems will not diverge arbitrarily far from each other due to varying state errors. Indeed, when the error is zero for all aggregated systems in a group, i.e. $x = \hat{x}$, \eqref{tempInterface} reduces to $\tilde{v}^\ell = (b/b_\ell) \cdot \hat{u}$. Thus, taking into account the difference between each $b_\ell$ and $b$, one can again use the bound on the error between $x$ and $\hat{x}$ to bound the deviation of $\tilde{v}^\ell$ from $\hat{u}$.}

\begin{table}[t]
	\caption{Partitioning of the 30 rooms into 3 groups.}
	\vspace{3.5mm}
	\label{partitioning}
	\centering
		\begin{tabular}{l | r r r}
			\toprule
			 				& Group 1 & Group 2 & Group 3 \\
			\midrule
			Pre-assignments 	& 1-6 & 11-18 & 21-27 \\
			Final partition	& 1-6 & 7-20 & 21-30 \\
			\bottomrule
		\end{tabular}
\end{table}

\subsection{Temperature Regulation}
We consider the task of regulating the temperature in a network of $L = 30$ rooms connected in a circle. The coupling matrix $\tilde{M} \in \R^{30 \times 30}$ is as shown in \eqref{exampleLaplacian}. We assume rooms 1-6, 11-18, and 21-27 are pre-assigned to 3 separate groups; the remaining rooms are assumed to be flexible with regard to temperature level, and are assigned to groups automatically via Optimization Problem 4. The pre-assignments and final partition are shown in Table \ref{partitioning}. The aggregate coupling matrix between the groups, obtained simultaneously with the final partition via Optimization Problem 4, is given by
\begin{equation}
\bar{M} = \begin{bmatrix}
-\sfrac{1}{3} & \sfrac{1}{6} & \sfrac{1}{6} \\
\sfrac{1}{14} & -\sfrac{1}{7} & \sfrac{1}{14} \\
\sfrac{1}{10} & \sfrac{1}{10} & -\sfrac{1}{5}
\end{bmatrix}. \nonumber
\end{equation}
One notes that this partition is \textit{not} equitable - indeed, with the pre-assignments shown in Table \ref{partitioning}, an equitable partition cannot be achieved. This is not problematic, however, since Theorem \ref{tempTheorem} relaxes the requirement of equitability of our partition, as long as we can find a matrix $Z \geq 0$ satisfying \eqref{relaxedCompositionalLMI}, where $Y = \bar{Y} \otimes I_p$ and $\mu_i = 1$, $i = 1, 2, 3$. Lemma \ref{ybarLemma} and Theorem \ref{tempTheorem} guarantee this is possible, however, since $\mathcal{N}(\tilde{M} + \tilde{M}^T)$ is spanned by $\mathbf{1}_L$ in this case, as $\tilde{M}$ is a Laplacian matrix. Thus, we solve Optimization Problem 3, with the additional constraint $\mu_i = 1$, $i = 1, 2, 3$ as mentioned, and obtain
\begin{equation}
Z = \begin{bmatrix}
2.0016 & -1.0490 & -0.9526 \\
-1.0490 & 1.9897 & -0.9407 \\
-0.9526 & -0.9407 & 1.8933
\end{bmatrix}. \nonumber
\end{equation}
Since we also relaxed the assumption of identical agents, the conduction coefficients $a_{\ell}$ and $b_{\ell}$ in our concrete model are permitted to vary between rooms. For each room, we select $a_{\ell}$ from a normal distribution with mean $0.005$ and standard deviation $0.0015$, and select $b_{\ell}$ from a normal distribution with mean $0.035$ and standard deviation $0.0075$. \addedTwo{Furthermore, since Theorem 1 allows us to aggregate subsystems with non-equal initial conditions, we select the initial temperature for each concrete room from a normal distribution with mean $18$ and standard deviation $0.15$. We then set
\begin{align}
\hat{\xi}_i(0) = \frac{1}{L_i} \left( \sum_{\ell \, \in \, \mbox{\small group}\ i} \xi^\ell(0) \right)
\end{align}
that is, the initial temperature of each aggregate room is equal to the average temperature of the aggregated rooms in its group. To demonstrate the robustness of our approach, we chose the standard deviation for the parameters and the initial conditions to be sufficiently large so that room temperatures within each group deviate visibly from each other during simulation (as seen in Figure \ref{temperatureSimulation}).}

We require the room temperature in the three areas of the building to increase to three separate temperature ranges in response to a signal which indicates, for example, that the building is currently occupied and must be adjusted to a more comfortable temperature. \added{This specification can be represented via, for example, a signal temporal logic (STL) formula \cite{donze2013signal,maler2004monitoring}. \addedTwo{Due to lack of space, we omit the details of the STL formula and refer the reader to \cite{smith2018hierarchical}, which includes two similar examples}. Although STL formulae are typically evaluated with respect to continuous time signals (see \cite{maler2004monitoring}, which considers dense-time real-valued signals), here we use the MPC approach from \cite{raman2017} which defines a semantics for STL over discrete time signals}. Since the approach in \cite{raman2017} utilizes mixed-integer programming, the computational burden of control synthesis of $\hat{\upsilon}$ is reduced significantly by using an aggregate model. The aggregate input is refined to a concrete input via the interface function \eqref{tempInterface} with $k_\ell = 2.5$ for $\ell = 1, \dots, L$. Simulation results are shown in Figure \ref{temperatureSimulation}.

%We discretize the continuous dynamics for the concrete and aggregate models and consider the STL formula
%\begin{equation} \label{stl1}
%\square( \varphi \wedge (S \implies \lozenge_{ [0, T] } \psi)) \nonumber
%\end{equation}
%\added{to be evaluated with respect to the time-step of $\tau = 1$ minute}. Here, $\varphi$ imposes the following temperature comfort bounds and input constraints
%\begin{align}
% \varphi := &(|\hat{\xi}_1 - 19| \ \leq 3) \wedge (|\hat{\xi}_2 - 19| \ \leq 3) \nonumber \\
% & (|\hat{\xi}_3 - 19| \ \leq 3) \wedge \ (0 \leq \hat{\upsilon}_1 \leq 5) \nonumber \\
% & (0 \leq \hat{\upsilon}_2 \leq 5) \wedge (0 \leq \hat{\upsilon}_3 \leq 5) 
%\varphi := \bigwedge_{i = 1}^3 \left[ (|\hat{\xi}_i - 19| \ \leq 3) \wedge (0 \leq \hat{\upsilon}_i \leq 5) \right] \nonumber \label{temperatureSafety}
%\end{align}
%\added{which are required to hold at each time-step}, and $\psi$ encodes a temperature target set
%\begin{equation} % \label{temperatureTarget}
%\psi := (18 \leq \hat{\xi}_1) \wedge (19 \leq \hat{\xi}_2) \wedge (20 \leq \hat{\xi}_3) \nonumber
%\end{equation}
%which must be reached within $T = 20$ time-steps after the signal $S$ is triggered. \added{In general, we note that by using MPC we are effectively imposing the state and input constraints on time-sampled versions of their respective signals.} 

\section{Conclusion} \label{conclusionSection}
In this paper we proposed to relax previous conditions required to construct an infinite abstraction for a non-stochastic dynamical system. We introduced a notion of practical simulation functions, which takes into account our relaxation and bounds the error between the concrete and abstract control systems. For a monolithic construction, we demonstrated that one can obtain a practical simulation function relating a linear control system to its abstraction, without requiring any geometric conditions to be satisfied. In the compositional case, we introduced a notion of practical storage functions, and showed how one can construct an abstraction and practical simulation function for an interconnected control system, without requiring a condition on the interconnection topology. In an application to aggregation, our theory enabled us to relax the assumption of identical agent models and equitability of the partition of the agents. We demonstrated this with a temperature regulation example, where the rooms in the building each have slightly varying dynamical models, and a non-equitable partition is used for aggregation.

%\section{Acknowledgements} ... \\
%\Majid{Many equations below do not require any number because you are not referring to them! This comment will apply to other parts of the paper as well!}

\appendix
\section{Appendix}
\subsection{Proof of Theorem \ref{linearTheorem}}
Let $\epsilon = x - P \hat{x}$ and note that we have the following bounds
\begin{align*}
\| h(x) - \hat{h}(\hat{x}) \|^2 = \epsilon^T C^T C \epsilon \leq \lambda_{\max} (C^T C) \|\epsilon\|^2, \quad \lambda_{\min}(U) \|\epsilon\|^2 \leq \epsilon^T U \epsilon = V(x, \hat{x}), % \label{outputproof2}
\end{align*}
for all $x$ and $\hat{x}$, since $\hat{C} = CP$. Thus, \eqref{output} holds with $\nu(s) = s^2 \lambda_{\min}(U) / \lambda_{\max}(C^T C)$, where $\nu \in \mathcal{K}_\infty$ since $U$ is positive definite.

Next, we apply the congruency transformation $\text{diag}(U, I)$ to \eqref{lmiconstr1}, yielding the equivalent condition
\begin{equation} \label{MKlmi}
\begin{bmatrix}
A_K^T U + U A_K + \alpha U & U W \\
W^T U & -\alpha I
\end{bmatrix} \leq 0. \nonumber
\end{equation}
where we have defined $A_K \triangleq A + BK$. Thus, for all $x$, $\hat{x}$, and $\hat{u}$ (determining $\epsilon$ and $d$), we have
\begin{align}
\begin{bmatrix}
\epsilon \\
d
\end{bmatrix}^T \begin{bmatrix}
A_K^T U + U A_K + \alpha U & U W \\
W^T U & -\alpha I
\end{bmatrix} \begin{bmatrix}
\epsilon \\
d
\end{bmatrix} \leq 0 \nonumber %\label{decayLMI}
\end{align} so that
\begin{align}
\nabla V\left(x, \hat{x}\right)^T\begin{bmatrix}
Ax+BK(x-P\hat x)+BQ\hat x+BR\hat u \\
\hat A\hat x+\hat B\hat u
\end{bmatrix} &= \epsilon^T \left[ A_K^T U + U A_K \right] \epsilon + d^T W^T U \epsilon + \epsilon^T U W d \nonumber \\
& \leq -\alpha \epsilon^T U \epsilon + \alpha d^T d \nonumber \\
& = -\alpha V(x, \hat{x}) + \alpha \| \hat{u} \|^2 + \alpha \| D \hat{x} \|^2 \nonumber %\label{decayINEQ}
\end{align}
which verifies that \eqref{iss} holds with $\eta(s) = \alpha s$, $\rho(s) = \alpha s^2$, and $\Delta(\hat{x}) = \alpha \| D \hat{x} \|^2$.

\subsection{Proof of Theorem \ref{compositionalityResultmod}}
Without modifications due to our relaxation, we can construct a $\mathcal{K}_\infty$ function $\nu$ satisfying \eqref{output} as in the proof of Theorem 4.2 given in \cite{ZA3}. Thus, we omit this portion of the proof and focus on showing that \eqref{iss} holds. We define the following error between the concrete and aggregate systems
\begin{equation}
\begin{bmatrix}
e_1 \\ \vdots \\ e_N
\end{bmatrix} :=
\begin{bmatrix}
h_{21}(x_1) - H_1 \hat{h}_{21}(\hat{x}_1) \\
\vdots \\
h_{2N}(x_N) - H_N \hat{h}_{2N}(\hat{x}_N)
\end{bmatrix}. \nonumber
\end{equation}
Then, from \eqref{interconnection} and \eqref{compositionalResidual}, it follows that
\begin{equation} \label{errorRelation}
W \begin{bmatrix} w_1 \\ \vdots \\ w_N \end{bmatrix}
- \hat{W} \begin{bmatrix} \hat{w}_1 \\ \vdots \\ \hat{w}_N \end{bmatrix} =
WM \begin{bmatrix} e_1 \\ \vdots \\ e_N \end{bmatrix} + Y
\begin{bmatrix}
\hat{h}_{21}(\hat{x}_1) \\
\vdots \\
\hat{h}_{2N}(\hat{x}_N)
\end{bmatrix}.
\end{equation}
Now, using the relation \eqref{errorRelation}, we obtain
\begin{align*}
& \begin{bmatrix}
W \begin{bmatrix} w_1 \\ \vdots \\ w_N \end{bmatrix}
- \hat{W} \begin{bmatrix} \hat{w}_1 \\ \vdots \\ \hat{w}_N \end{bmatrix} \\
h_{21}(x_1) - H_1 \hat{h}_{21}(\hat{x}_1) \\
\vdots \\
h_{2N}(x_N) - H_N \hat{h}_{2N}(\hat{x}_N)
\end{bmatrix}^T X(\mu_1 X_1, \dots, \mu_N X_N) \begin{bmatrix}
W \begin{bmatrix} w_1 \\ \vdots \\ w_N \end{bmatrix}
- \hat{W} \begin{bmatrix} \hat{w}_1 \\ \vdots \\ \hat{w}_N \end{bmatrix} \\
h_{21}(x_1) - H_1 \hat{h}_{21}(\hat{x}_1) \\
\vdots \\
h_{2N}(x_N) - H_N \hat{h}_{2N}(\hat{x}_N)
\end{bmatrix} \\
= &\begin{bmatrix}
\hat{h}_{21}(\hat{x}_1) \\
\vdots \\
\hat{h}_{2N}(\hat{x}_N) \\
e_1 \\
\vdots \\
e_N \\
\end{bmatrix}^T \begin{bmatrix}
Y & W M \\
0 & I_{\tilde{q}}
\end{bmatrix}^T
X
\begin{bmatrix}
Y & W M \\
0 & I_{\tilde{q}}
\end{bmatrix} \begin{bmatrix}
\hat{h}_{21}(\hat{x}_1) \\
\vdots \\
\hat{h}_{2N}(\hat{x}_N) \\
e_1 \\
\vdots \\
e_N \\
\end{bmatrix} \leq \begin{bmatrix}
\hat{h}_{21}(\hat{x}_1) \\
\vdots \\
\hat{h}_{2N}(\hat{x}_N)
\end{bmatrix}^T Z
\begin{bmatrix}
\hat{h}_{21}(\hat{x}_1) \\
\vdots \\
\hat{h}_{2N}(\hat{x}_N)
\end{bmatrix}
\end{align*}
where the inequality follows from the fact that $Z$ and $\mu_1, \dots, \mu_N$ satisfy \eqref{relaxedCompositionalLMI}. Using this bound, the proof of Theorem 4.2 given in \cite{ZA3} can be easily modified to show that \eqref{iss} holds for an appropriate choice of $\eta \in \mathcal{K}_\infty$, $\rho \in \mathcal{K} \cup \{0\}$, and with $\Delta(\hat{x})$ as defined in \eqref{nonequitableDeviation}. Therefore, we conclude that $V$ in \eqref{compositionalSimFun} is a practical simulation function from $\hat{\Sigma}$ to $\Sigma$.

\subsection{Proof of Lemma \ref{ybarLemma}}
We note that $P$ has the form $P=\text{diag}(\mathbf{1}_{L_1}, \cdots, \mathbf{1}_{L_N})$
%\begin{equation}
%P = \begin{bmatrix}
%\mathbf{1}_{L_1} & & \\
%& \ddots & \\
%& & \mathbf{1}_{L_N}
%\end{bmatrix} \nonumber
%\end{equation}
upon a permutation. Therefore,
\begin{equation}
P \bar{M} = \begin{bmatrix}
\bar{m}_{11} \mathbf{1}_{L_1} & \dots & \bar{m}_{1N} \mathbf{1}_{L_1} \\
\bar{m}_{21} \mathbf{1}_{L_2} & \dots & \bar{m}_{2N} \mathbf{1}_{L_2} \\
\vdots & & \vdots \\
\bar{m}_{N1} \mathbf{1}_{L_N} & \dots & \bar{m}_{NN} \mathbf{1}_{L_N}
\end{bmatrix} \nonumber
\end{equation}
where the $\bar{m}_{ij} \in \R$ denote entries of $\bar{M}$. Let
\begin{equation}
\begin{bmatrix}
v_{11} & \dots & v_{1N} \\
v_{21} & \dots & v_{2N} \\
\vdots & & \vdots \\
v_{N1} & & v_{NN}
\end{bmatrix} := \tilde{M} P, \quad v_{ij} \in \R^{L_i}. \nonumber
\end{equation}
Then, we see that
\begin{equation}
\tilde{M} P - P \bar{M} = \begin{bmatrix}
v_{11} - \bar{m}_{11} \mathbf{1}_{L_1} & \dots & v_{1N} - \bar{m}_{1N} \mathbf{1}_{L_1} \\
\vdots & & \vdots \\
v_{N1} - \bar{m}_{N1} \mathbf{1}_{L_N} & \dots & v_{NN} - \bar{m}_{NN} \mathbf{1}_{L_N}
\end{bmatrix} \nonumber
\end{equation}
and
\begin{equation}
\|\tilde{M} P - P \bar{M} \|_F = \left\| \begin{bmatrix}
v_{11} - \bar{m}_{11} \mathbf{1}_{L_1} \\
\vdots \\
v_{N1} - \bar{m}_{N1} \mathbf{1}_{L_N} \\
%v_{12} - \bar{m}_{12}\mathbf{1}_{L_1} \\
\vdots \\
%v_{N2} - \bar{m}_{N2} \mathbf{1}_{L_N} \\
%\vdots \\
v_{1N} - \bar{m}_{1N} \mathbf{1}_{L_1} \\
\vdots \\
v_{NN} - \bar{m}_{NN} \mathbf{1}_{L_N}
\end{bmatrix} \right\|. \nonumber
\end{equation}
Minimization of the latter Euclidean norm over $\bar{M}$ can be decomposed into the independent problems
\begin{equation}
\min_{\bar{m}_{ij}} \|v_{ij} - \bar{m}_{ij} \mathbf{1}_{L_i}\|, \quad i,j = 1, \dots, N. \nonumber
\end{equation}
Since $\|v_{ij} - \bar{m}_{ij} \mathbf{1}_{L_i}\|^2 = (v_{ij} - \bar{m}_{ij} \mathbf{1}_{L_i})^T (v_{ij} - \bar{m}_{ij} \mathbf{1}_{L_i}) = v_{ij}^T v_{ij} - 2 \bar{m}_{ij} \mathbf{1}_{L_i}^T v_{ij} + \bar{m}_{ij}^2 L_i$, the minimizer is $\bar{m}_{ij}^* = (1/L_i) \mathbf{1}_{L_i}^T v_{ij}$.

We now verify the claim of Lemma \ref{ybarLemma}; we have
\begin{equation}
\mathbf{1}_L^T P = \mathbf{1}_L^T
\begin{bmatrix}
\mathbf{1}_{L_1} & & \\
& \ddots & \\
& & \mathbf{1}_{L_N}
\end{bmatrix} = \begin{bmatrix} L_1 & \dots & L_N \end{bmatrix} \nonumber
\end{equation}
thus,
\begin{align*}
\mathbf{1}_L^T P \bar{M} &= \begin{bmatrix}
L_1 & \dots & L_N
\end{bmatrix} \bar{M} = \begin{bmatrix}
\sum_{i=1}^N \bar{m}_{i1} L_i & \dots & \sum_{i=1}^N \bar{m}_{iN} L_i \end{bmatrix}.
\end{align*}
Since the optimal values for $\bar{m}_{ij}$ give
\begin{equation}
\sum_{i=1}^N \bar{m}^*_{ij} L_i = \sum_{i=1}^N \mathbf{1}^T_{L_i} v_{ij} = \mathbf{1}_L^T 
\begin{bmatrix}
v_{1j} \\
\vdots \\
v_{Nj}
\end{bmatrix}, \nonumber
\end{equation}
we get
\begin{equation}
\mathbf{1}_L^T P \bar{M} = \mathbf{1}_L^T \begin{bmatrix}
v_{11} & \dots & v_{1N} \\
\vdots & & \vdots \\
v_{N1} & & v_{NN}
\end{bmatrix} = \mathbf{1}_L^T \tilde{M} P \nonumber
\end{equation}
and therefore $\mathbf{1}_L^T (P \bar{M} - \tilde{M} P) = \mathbf{1}_L^T \bar{Y} = 0$.

\subsection{Proof of Proposition \ref{Prop2}}
If we let
$
\tilde{V}^\ell=\frac{1}{2}(x-\hat{x})^TQ_\ell(x-\hat{x}),
$ then (\ref{lowerbound}) holds with $\varsigma(x)=x$, $\tilde{\nu}^\ell(s)=\frac{1}{2}\lambda_{\rm min}(Q_\ell)s^2$, and (\ref{derivative}) becomes
\begin{align} \nonumber
& (x-\hat{x})^TQ_\ell (\alpha_\ell(x)-\alpha(\hat{x}))+(x-\hat{x})^TQ_\ell\left(\beta_\ell(x)\tilde{v}^\ell(x,\hat{x},\hat{u})-\beta(\hat{x})\hat{u}\right)+(x-\hat{x})^TQ_\ell B(w-\hat{w}) \\
& \le (x-\hat{x})^TQ_\ell (\alpha_\ell(x)-\alpha(\hat{x}))+\vartheta_\ell \|x-\hat{x}\|^2 +\tilde{\rho}^\ell(\|\hat{u}\|) +(\sigma(x)-\sigma(\hat{x}))^T(w-\hat{w}), \label{intm1}
\end{align}
where the inequality follows from (\ref{Cond2}) and (\ref{Cond4}), combined with $\sigma(x)=Cx$ from (\ref{specialclass}). We rewrite the first term on the right hand side of (\ref{intm1}) as 
\begin{equation}
(x-\hat{x})^TQ_\ell (\alpha_\ell(x)-\alpha(\hat{x})) =(x-\hat{x})^TQ_\ell (\alpha_\ell(x)-\alpha_\ell(\hat{x})) +(x-\hat{x})^TQ_\ell (\alpha_\ell(\hat{x})-\alpha(\hat{x})). \label{intm2}
\end{equation}
It follows from (\ref{Cond1}) that 
\begin{equation}\label{intm3}
 (x-\hat{x})^TQ_\ell (\alpha_\ell(x)-\alpha_\ell(\hat{x})) \le \lambda_\ell \|x-\hat{x}\|^2.
\end{equation}
 To see this, define the function $\Omega(t)=\alpha_\ell(\hat{x}+t(x-\hat{x}))$ and note
 \begin{equation}\label{calc}
 (x-\hat{x})^TQ_\ell\int_0^1 \Omega'(t)dt
 \end{equation}
 is equal to the left hand side of (\ref{intm3}) by the fundamental theorem of calculus. From the chain rule, (\ref{calc}) equals 
  \begin{equation}\label{calc2}
 (x-\hat{x})^TQ_\ell\int_0^1 J(\hat{x}+t(x-\hat{x}))dt(x-\hat{x})
 \end{equation}
 where $J$ is the Jacobian of $\alpha_\ell$.  Rewriting (\ref{calc2}) as
   \begin{equation}
\frac{1}{2} (x-\hat{x})^T \left( \int_0^1 \left( Q_\ell J + J^TQ_\ell \right) dt \right) (x-\hat{x}), \nonumber
 \end{equation}
 we see from (\ref{Cond1}) that the integrand is bounded above by $2\lambda_\ell I_n$, which confirms (\ref{intm3}).
Next, we note that
 \begin{equation}
(x-\hat{x})^T Q_\ell (\alpha_\ell(\hat{x})-\alpha(\hat{x})) \le \kappa \|x-\hat{x}\|^2 +\frac{1}{4\kappa}\|Q_\ell (\alpha_\ell(\hat{x})-\alpha(\hat{x})) \|^2 \label{intm4}
\end{equation}
for any choice of $\kappa>0$, which follows from Young's inequality \cite{young1912classes}. % as can be shown with a completion of squares argument \Majid{Do you mean Young's inequality?}.
Then, from (\ref{intm2}), (\ref{intm3}) and (\ref{intm4}), an upper bound on (\ref{intm1}) is
\begin{equation}
(\lambda_\ell +\vartheta_\ell+\kappa) \|x-\hat{x}\|^2 +\frac{1}{4\kappa}\|Q_\ell (\alpha_\ell(\hat{x})-\alpha(\hat{x})) \|^2+\tilde{\rho}^\ell(\|\hat{u}\|)+(\sigma(x)-\sigma(\hat{x}))^T(w-\hat{w}).\label{intm5}
\end{equation}
We select $\kappa=|\lambda_\ell +\vartheta_\ell |-\varepsilon_\ell$, which is positive since $\varepsilon_\ell \in (0, |\lambda_\ell+\vartheta_\ell|)$, and note that (\ref{intm5}) becomes
\begin{equation} -\varepsilon_\ell  \|x-\hat{x}\|^2+\frac{1}{4(|\lambda_\ell +\vartheta_\ell |-\varepsilon_\ell)}\|Q_\ell (\alpha_\ell(\hat{x})-\alpha(\hat{x})) \|^2+\tilde{\rho}^\ell(\|\hat{u}\|)+(\sigma(x)-\sigma(\hat{x}))^T(w-\hat{w}). \label{intm6}
 \end{equation}
Substituting the inequality $\varepsilon_\ell  \|x-\hat{x}\|^2\ge \frac{2\varepsilon_\ell}{\lambda_{\rm max}(Q_\ell)}\tilde{V}^\ell=\tilde{\eta}^\ell(\tilde{V}^\ell)$ in (\ref{intm6}), we obtain (\ref{agent_dissip}) with the terms defined in (\ref{passive}).

\subsection{Proof of Theorem \ref{tempTheorem}}
We have shown the equitability criterion \eqref{equit} is identical to condition \eqref{equitabilityconstr} of Theorem \ref{compositionalityResult}; also, that if we select $\mu_i = 1$, $i = 1, \dots, N$, then \eqref{interconnectionCondition} implies condition \eqref{equitableLMI} of Theorem \ref{compositionalityResult} holds. Thus, if we use an equitable partition for aggregation and \eqref{interconnectionCondition} holds, then both conditions of Theorem \eqref{compositionalityResult} also hold so that \eqref{compositionalSimFun} is indeed a simulation function from $\hat{\Sigma}$ to $\Sigma$, with $V_i(x_i, \hat{x}_i)$ as in \eqref{storage}, and where $\mu_i = 1$, $i = 1,\dots,N$. It follows that relaxing the equitability condition as in \eqref{equitRelaxation} is identical to the relaxation \eqref{compositionalResidual} given in Theorem \ref{compositionalityResultmod}. Thus, in this case one must choose a matrix $Z = Z^T \geq 0$ satisfying \eqref{relaxedCompositionalLMI}, with $Y = \hat{Y} \otimes I_p$ and $\mu_i = 1$ for $i = 1,\dots,N$.

To show that $\mathcal{N}(\tilde{M} + \tilde{M}^T) \subseteq \mathcal{N}(\bar{Y}^T)$ is a necessary and sufficient condition for such a $Z$ to exist, we prove the following fact. Let $B \in \R^{m \times n}$ be an arbitrary matrix and $C \in \R^{n \times n}$ be such that $C = C^T \leq 0$. Then, there exists a matrix $A \in \R^{n \times n}$ such that $A = A^T \geq 0$ and
\begin{equation}
\begin{bmatrix}
-A & B \\
B^T & C
\end{bmatrix} \leq 0 \label{simpleCondition}
\end{equation}
if and only if $\mathcal{N}(C) \subseteq \mathcal{N}(B)$. To see the necessity, suppose there exists a vector $y \in \R^{n}$ such that $y \in \mathcal{N}(C)$ but $y \notin \mathcal{N}(B)$. Then, for any $x \in \R^m$, we have
\begin{align}
\begin{bmatrix}
x \\ y
\end{bmatrix}^T
\begin{bmatrix}
-A & B \\
B^T & C
\end{bmatrix}
\begin{bmatrix}
x \\ y
\end{bmatrix} = -x^T A x + 2 x^T B y. \label{necessityProof}
\end{align}
Let $x = \theta B y$, where $\theta \in \R_{>0}$. Then, a lower bound for \eqref{necessityProof} is
\begin{equation}
(2\theta - \theta^2 \lambda_{\text{max}}(A))\|By\|^2 \nonumber
\end{equation}
which is positive for any choice of $\theta \in (0, 2/\lambda_{\text{max}}(A))$. Thus, for any $A = A^T \geq 0$, condition \eqref{simpleCondition} does not hold. For the sufficiency, suppose $\mathcal{N}(C) \subseteq \mathcal{N}(B)$, and let $\phi > 0$ be the smallest nonzero eigenvalue of $-C$ (if $-C$ has no nonzero eigenvalues, then $C$ is the zero matrix and the proof follows trivially). We select $A = -(1/\phi) B B^T$, and note that
\begin{equation}
\begin{bmatrix}
x \\ y
\end{bmatrix}^T
\begin{bmatrix}
-(1/\phi) B B^T & B \\
B^T & C
\end{bmatrix}
\begin{bmatrix}
x \\ y
\end{bmatrix} = -(1/\phi) x^T B B^T x + 2 x^T B y + y^T C y. \label{sufficiencyProof} 
\end{equation}
Next, we decompose $y$ as $y = y_1 + y_2$, where $y_1 \in \mathcal{N}(C)$ and $y_1^T y_2 = 0$. We note that, by assumption, \eqref{sufficiencyProof} becomes
\begin{align}
-(1/\phi) x^T B B^T x + 2 x^T B y_2 + y_2^T C y_2 &\leq -(1/\phi) x^T B B^T x + 2 x^T B y_2 - \phi \|y_2\|^2 \nonumber \\
&= -(1/\phi) \|z\|^2 + 2 z^T y_2 - \phi \|y_2\|^2 \label{sufficiencyFinal}
\end{align}
where the second step follows since $y_2 \notin \mathcal{N} (C)$, and the third step results from the definition $z := B^T x$. Finally, using Young's inequality \cite{young1912classes} as
\begin{equation}
z^T y_2 \leq \frac{1}{2 \phi} \| z \|^2 + \frac{\phi}{2} \| y_2 \|^2 \nonumber
\end{equation}
one can see that \eqref{sufficiencyFinal} is bounded above by zero.

We can then recover the null space condition of Theorem 5 as follows. Using the notation in \eqref{relaxedCompositionalLMI}, we note that
\begin{equation}
Q(Z,1,\dots,1)
= \begin{bmatrix}
-Z & \frac{1}{2} \bar{Y}^T \otimes I_p \\
\frac{1}{2} \bar{Y} \otimes I_p & \frac{1}{2} (\tilde{M} + \tilde{M}^T) \otimes I_p
\end{bmatrix} \nonumber
\end{equation}
which can be mapped to the matrix in \eqref{simpleCondition} by taking $A = Z$, $B = \frac{1}{2}\bar{Y}^T \otimes I_p$ and $C = \frac{1}{2}(\tilde{M} + \tilde{M}^T) \otimes I_p$. Thus, a necessary and sufficient condition for the existence of a $Z = Z^T \geq 0$ such that $Q(Z,1,\dots,1) \leq 0$ is $\mathcal{N}(\frac{1}{2}(\tilde{M} + \tilde{M}^T) \otimes I_p) \subseteq \mathcal{N}(\frac{1}{2}\bar{Y}^T \otimes I_p)$, which is equivalent to $\mathcal{N}(\tilde{M} + \tilde{M}^T) \subseteq \mathcal{N}(\bar{Y}^T)$.

% Suppose there exists a vector $y \in \R^{L}$ such that $y \in \mathcal{N}(\tilde{M} + \tilde{M}^T)$ but $y \notin \mathcal{N}(\bar{Y}^T)$. Then, using the notation in \eqref{relaxedCompositionalLMI}, we have for any $x \in \R^{Np}$
% \begin{align}
% &\begin{bmatrix} x \\ y \otimes \mathbf{1}_p \end{bmatrix}^T
% Q(Z,1,\dots,1) \begin{bmatrix} x \\ y \otimes \mathbf{1}_p \end{bmatrix} \nonumber \\
% =& \begin{bmatrix} x \\ y \otimes \mathbf{1}_p \end{bmatrix}^T
% \begin{bmatrix}
% -Z & \frac{1}{2} \bar{Y}^T \otimes I_p \\
% \frac{1}{2} \bar{Y} \otimes I_p & \frac{1}{2} (\tilde{M} + \tilde{M}^T) \otimes I_p
% \end{bmatrix} \begin{bmatrix} x \\ y \otimes \mathbf{1}_p \end{bmatrix} \nonumber \\
% =& -x^T Z x + x^T \left[ (\bar{Y}^T y) \otimes \mathbf{1}_p \right].
% \end{align}
% Let $x = \alpha (\bar{Y}^T y) \otimes \mathbf{1}_p$, where $\alpha \in \R$. Then, a lower bound for \eqref{necessityProof} is
% \begin{align}
% (\alpha - \alpha^2 \lambda_{\text{max}}(Z)) \| (\bar{Y}^T y) \otimes \mathbf{1}_p \|^2
% \end{align}
% which is positive for any choice of $\alpha \in (0,1/ \lambda_{\text{max}}(Z))$. Thus, for any $Z = Z^T \geq 0$, condition \eqref{relaxedCompositionalLMI} is violated. To see the sufficiency, 

% \Majid{Please unify references! In some you use full name of authors and in some you use initial of first name.}

\bibliographystyle{plain}
\bibliography{bibfile}

\end{document}